\newif\ifpdf \ifx\pdfoutput\undefined \pdffalse \else \pdfoutput=1
\pdftrue \fi

 \newenvironment{fig}[2] { \begin{center}
 \begin{figure}
 \def\myfigcaptionlar{#1} \def\myfigcaptionvar{{#2}}}
 {\caption{\myfigcaptionvar}\label{\myfigcaptionlar}
 \end{figure}\end{center}}
 
 \newcommand{\psfig}[3] {\begin{fig}{#1}{#2}\if 0#3 \else
 \begin{center}\includegraphics{#3}\end{center}\fi
 \end{fig} }
 
 \newcommand{\texfig}[3] {\begin{fig}{#1}{#2} \if 0#3 \else
 \begin{center}
 \input{#3}
 \end{center}
 \fi
 \end{fig}
 \if 0#3 \fi  }
 
 \newcommand{\citez}[1]{\cite{#1}} 

\documentclass[openany,oneside,fleqn,11pt]{article}
\usepackage{amsmath,amssymb,amscd,amsthm} 
\usepackage{url} 
\usepackage[numbers]{natbib}
 \usepackage{float}
\ifpdf \usepackage[pdftex]{graphicx}
\else 
 \usepackage{graphicx} 
\fi
\usepackage[bookmarks=true]{hyperref} %this one must be the last one
\newcommand{\hyref}[1]{\href{#1}{\tt{#1}}}
\newcommand{\hyreff}[2]{\href{#1}{\tt{#2}}}
 
\begin{document}
%\DeclareGraphicsExtensions{.pdf} 
\DeclareGraphicsExtensions{.mps}

%%%%%%%%%%%%%%% Textual Macros %%%%%%%%%%%%%%%%%%%%%%%%%%%%%%%%%%%%%%%%%%%%
\newcommand{\Ito}{It\^o }

\newcommand{\Mathematica}{{\it
 Mathematica }}
\newcommand{\MathReader}{{\it
 MathReader }}

\newcommand{\abbrev}[1]{#1. }
\newcommand{\ie}{\abbrev{i.e}}
\newcommand{\pg}{\abbrev{p}}
\newcommand{\eg}{\abbrev{eg}}
\newcommand{\etc}{\abbrev{etc}}
\newcommand{\vs}{\abbrev{vs}}
\newcommand{\propan}{\abbrev{Propos}}
\newcommand{\theopan}{Theorem}

%%%%%%%%%%%%%%% Propositions %%%%%%%%%%%%%%%%%%%%%%%%%%%%%%%%%%%%%%%%%%%%%%%%%%
\theoremstyle{theorem} 
\newtheorem{prop}{Proposition}
\newtheorem{lemma}[prop]{Lemma} 
\newtheorem{theo}[prop]{Theorem}
\newtheorem{coroll}[prop]{Corollary}
\newcommand{\refprop}[1]{Proposition~\ref{#1}}
\newcommand{\refppropag}[1]{\propan~\ref{#1} \pg~\pageref{#1}}
\newcommand{\reftpropag}[1]{\theopan~\ref{#1} \pg~\pageref{#1}}
\newcommand{\refcorol}[1]{Corollary~\ref{#1} \pg~\pageref{#1}}
\theoremstyle{definition} 
\newtheorem{adefn}[prop]{Definition}
\newenvironment{defn}{\begin{adefn}}{\hfill$\blacksquare$\end{adefn}}
\theoremstyle{remark}
\newtheorem{arem}[prop]{Remark}
\newtheorem{aexample}[prop]{Example}
\newcommand{\defbf}[1]{{\bf #1}}
\newenvironment{solution}{\begin{proof}[Solution]}{\end{proof}}
\newenvironment{rem}{\begin{arem}}{\hfill$\triangledown$\end{arem}}
\newenvironment{example}{\begin{aexample}}{\hfill$\vartriangle$\end{aexample}}

\newcommand{\algitem}[1]{\item #1}

\floatstyle{boxed} \newfloat{Algorithm}{H}{lop}
\newenvironment{alg}[5]{ \def\myalgcaptionvar{{\bf #2}$\quad[$#1$]$}
 \if 0#5 \def\myalgcaptionvarl{\relax} \else
 \def\myalgcaptionvarl{\noindent{#5}} \fi
 \begin{Algorithm}
 \noindent{\bf input: }#3
 \noindent{\bf output: }#4
 \begin{itemize}
 } {
 \end{itemize}
 \myalgcaptionvarl
% \hrule
 \caption{\myalgcaptionvar}
 \end{Algorithm}
}
 
\newenvironment{romanenu}
{\begin{enumerate}\renewcommand{\labelenumi}{(\roman{enumi})}}
 {\end{enumerate}}

%%%%%%%%%%%%%%% Notation %%%%%%%%%%%%%%%%%%%%%%%%%%%%%%%%%%%%%%%%%%%%%%%%%%
% common ==========================
\newcommand{\rfor}[1]{\quad\text{for}\quad#1}
\newcommand{\firstfor}{\qquad&&\text{for}\quad}
\newcommand{\nextfor}{&&\text{for}\quad}

\newcommand{\alterbar}{\vec}

\newcommand{\origi}[1]{{#1}_0} \newcommand{\define}{:=}
\newcommand{\definer}{=:} \newcommand{\compose}{\circ}
\newcommand{\frahalf}{{\frac{1}{2}}}

% iteration =======================
\newcommand{\iva}[1]{0 \le i \le #1} \newcommand{\jva}[1]{0 \le j \le
 #1}
\newcommand{\ssum}[3]{{\sum_{#1 = #2}^{#3}}}
\newcommand{\ssumi}[1]{{\ssum{i}{1}{#1}}}

% sets and spaces =================
\newcommand{\RR}{{\mathbb R}} \newcommand{\NN}{{\mathbb N}}
\newcommand{\ZZ}{{\mathbb Z}} \newcommand{\CC}{{\mathbb C}}

\newcommand{\RntoRn}{{\colon\RR^n \to \RR^n}}
\newcommand{\RntoR}{{\colon\RR^n \to \RR}}
\newcommand{\RtoRn}{{\colon\RR \to \RR^n}}
\newcommand{\RtoR}{{\colon\RR \to \RR}}

% systems =========================
\newcommand{\dsys}[3]{\left(#1,#2,U,\origi{#3}\right)}
\newcommand{\ssys}[4]{\left(#1,#2,#3,U,\origi{#4}\right)}
\newcommand{\osys}[4]{\left(#1,#2,#3,U,\origi{#4}\right)}
\newcommand{\dsysfgx}{\dsys{f}{g}{x}}
\newcommand{\dosysfghx}{\osys{f}{g}{h}{x}}
\newcommand{\ssysfgx}{\ssys{f}{g}{\sigma}{x}}

\newcommand{\bigdsys}[3]{\left(#1\, , \, #2\, , \, U\, , \,
 \origi{#3}\right)} \newcommand{\bigssys}[4]{\left(#1\, , \, #2\, ,
 \, #3\, , \, U\, , \, \origi{#4}\right)}

\newcommand{\clssys}{{\mathbb X}} \newcommand{\clssysdet}{\clssys_D}
\newcommand{\clssysito}{\clssys_I}
\newcommand{\clssysstrat}{\clssys_S}

\newcommand{\clsys}{\clssys(n,1,1)}
\newcommand{\clsysdet}{\clssysdet(n,1)}
\newcommand{\clsysito}{\clssysito(n,1,1)}
\newcommand{\clsysstrat}{\clssysstrat(n,1,1)}

\newcommand{\clmisys}{\clssys(n,m,k)}
\newcommand{\clmisysdet}{\clssysdet(n,m)}
\newcommand{\clmisysito}{\clssysito(n,m,k)}
\newcommand{\clmisysstrat}{\clssysstrat(n,m,k)}

\newcommand{\dsdef}{$\Theta_D = \dsysfgx \in \clsysdet$ }
\newcommand{\isdef}{$\Theta_I = \ssysfgx \in \clsysito$ }
\newcommand{\ssdef}{$\Theta_S = \ssysfgx \in \clsysstrat$ }

\newcommand{\dsdefmiuu}{$\Theta = \ssysfgx \in \clmisys$ }
\newcommand{\dsdefmi}{$\Theta_D = \dsysfgx \in \clmisysdet$ }
\newcommand{\isdefmi}{$\Theta_I = \ssysfgx \in \clmisysito$ }
\newcommand{\ssdefmi}{$\Theta_S = \ssysfgx \in \clmisysstrat$ }

\newcommand{\smooth}{C^{\infty}} \newcommand{\analytic}{C^{\omega(x)}}
\newcommand{\kdifferentiable}{C^{k}}

% operations

% differentiation ============
\newcommand{\parby}[2]{\frac{\partial #1}{\partial #2}}
\newcommand{\parbyx}[1]{\parby{#1}{x}}
\newcommand{\parbyz}[1]{\parby{#1}{z}}

\newcommand{\parbysec}[2]{\frac{\partial^2 #1}{\partial #2^2}}

\newcommand{\parbyxx}[3]{\frac{\partial^2 #1}{\partial #2 \partial
 #3}} \newcommand{\parbyxxx}[4]{\frac{\partial^3 #1}{\partial #2
 \partial #3 \partial #4}}
\newcommand{\parbyxxxx}[5]{\frac{\partial^4 #1}{\partial #2 \partial
 #3 \partial #4 \partial #5 }}

\newcommand{\lie}[2]{{\mathcal L}_{#1} {#2}}
\newcommand{\multilie}[3]{{\mathcal L}_{#1}^{#2} {#3}}
\newcommand{\biglie}[2]{\langle d{#2},{#1}\rangle}

\newcommand{\ad}[3]{\operatorname{ad}_{{#1}}^{{#2}} {{#3}}}
\newcommand{\adfg}[1]{\ad{f}{{#1}}{g}}
\newcommand{\lighk}[1]{\lie{g}{\multilie{f}{#1}{h}}}
\newcommand{\distrofg}[1]{\left\{\ad{f}{i}{g},\,\iva{{#1}}\right\}}
\newcommand{\distroabafg}[1]{\left\{\ad{{\alterbar
 f}}{i}{g},\,\iva{{#1}}\right\}}

\newcommand{\compoit}{\compose T^{-1} (z)}
\newcommand{\travf}[1]{\parbyx{T} #1 \compoit}
\newcommand{\trav}[1]{\parbyx{T} #1}

% operators ==================
\newcommand{\flow}[1]{Fl^{#1}}
\newcommand{\trace}{{\operatorname{trace}}}
\newcommand{\grad}{{\operatorname{grad}}}
\newcommand{\kernel}{{\operatorname{kernel}}}
\newcommand{\annihilator}{{\operatorname{annihilator}}}
\newcommand{\rank}{{\operatorname{rank}}}
\newcommand{\abs}{{\operatorname{abs}}}
\newcommand{\sgn}{{\operatorname{sgn}}}
\newcommand{\sspan}{{\operatorname{span}}} \newcommand{\ito}[2]{P_#1
 #2} \newcommand{\corr}[2]{\operatorname{corr}_#1(#2)}
\newcommand{\Corr}[2]{\operatorname{Corr}_#1 #2}

\newcommand{\corrxx}[2]{{\frahalf \parby{#1}{#2} #1}}
\newcommand{\corrx}{\corrxx{\sigma}{x}}

\newcommand{\itox}{{\frahalf \sigma^2 \parbysec{T}{x} \compose T^{-1}
 (z)}} \newcommand{\itoxs}{{\frahalf \sigma^2 \parbysec{T}{x} }}

% coordinate transformation ======
\newcommand{\tantra}[1]{#1_\ast} \newcommand{\cotantra}[1]{#1^\ast}

\newcommand{\sct}[1]{{\mathcal T}_{#1}} \newcommand{\sctt}{\sct{T}}

\newcommand{\scti}[1]{\sct{#1}^{I}}
\newcommand{\scts}[1]{\sct{#1}^{S}}

\newcommand{\feedback}[2]{{\mathcal F}_{#1,#2}}
\newcommand{\feedbackab}{\feedback{{\alpha}}{{\beta}}}

\newcommand{\combinedz}{{\mathcal J}}
\newcommand{\combined}[3]{\combinedz_{#1,#2,#3}}
\newcommand{\combinedtab}{\combined{T}{{\alpha}}{{\beta}}}

% stochastic
\newcommand{\expected}[1]{{\cal E}\left\{#1\right\}}
\newcommand{\var}[1]{\operatorname{var}\left\{#1\right\}}
\newcommand{\cov}[1]{\operatorname{cov}\left\{#1\right\}}

\newcommand{\wrong}{^\#}

%%%%%%%%%%%%%%%%%%%%%%%%%%%%%%%%%%%%%%%%%%%%%%%%%%%%%%%%%%%%
\newcommand{\tpDate}{2002/10/09}
\newcommand{\tpAuthor}{Ladislav Sl\'ade\v cek}
\newcommand{\tpeAddress}{\v R\'\i{}kovice 18, CZ 751 18, Czech Republic}

\newcommand{\tpEmail}{\hyreff{mailto:lsla@post.cz}{lsla@post.cz}}
\newcommand{\tpeTitle}{Exact Feedback Linearization of Stochastic
Control Systems}
% #3
\newcommand{\tpeAbstractN}{Abstract}
% #4
\newcommand{\tpeAbstract}{This paper studies exact linearization
 methods for stochastic SISO affine controlled dynamical systems. The
 systems are defined as vectorfield triplets in Euclidean space. The
 goal is to find, for a given nonlinear stochastic system, a
 combination of invertible transformations which transform the system
 into a controllable linear form. Of course, for most nonlinear
 systems  such transformation does not exist.

  We are focused on linearization by state coordinate transformation
 combined with feedback.  The difference between \Ito and Stratonovich
 systems is emphasized. Moreover, we define three types of linearity
 of stochastic systems --- $g$-linearity, $\sigma$-linearity, and
 $g\sigma$-linearity.
 
 Six variants of the stochastic exact linearization problem are
 studied. The most useful problem --- the \Ito-~$g\sigma$
 linearization is solved using the correcting term, which proved to be
 a very useful tool for \Ito systems.  The results are illustrated on
 a numerical example solved with help of symbolic algebra.  }
% #5
\newcommand{\tpeKeywordsN}{Keywords}
% #6
\newcommand{\tpeKeywords}{ exact linearization, feedback
  linearization, nonlinear dynamical system, \Ito integral,
  Stratonovich integral, correcting term 
{\bf MCS classification:} 93B18, 93E03}

% ============================================================
\title{\tpeTitle} \author{\tpAuthor}

\newcommand{\tp}[6]{ \vfill
 \begin{titlepage}
 \begin{center} 
 {~\\} \vspace{25mm}
 {\Large #1}\\
 \vspace{1cm}
 {\large \tpAuthor}\\
 {\tpEmail}\\
 {#2}\\
 
 \vspace{0.5cm}
 {\large \tpDate}\\
 \end{center}
 \begin{quote}\small
 \vspace{15mm}
 {\sc\bf #3:} #4\\
 {\sc\bf #5:} #6\\
 \end{quote}
 \end{titlepage}
 }
 
 \newcommand{\tpe}{ \tp {\tpeTitle} {\tpeAddress} {\tpeAbstractN}
 {\tpeAbstract} {\tpeKeywordsN} {\tpeKeywords} }

%%%%%%%%%%%%%%%%%%%%%%%%%%%%%%%%%%%%%%%%%%%%%%%%%%%%%%%%%%%%
\tpe 

\tableofcontents
%%%%%%%%%%%%%%% Section %%%%%%%%%%%%%%%%%%%%%%%%%%%%%%%%%%%%%%%%%%%%%

\section{Introduction}
\pagenumbering{arabic}\count1=0
\label{sec:intro}

The theory of exact linearization of deterministic dynamical systems
has been thoroughly studied since seventies. This paper attempts to
apply some of the results to the stochastic area. We emphasize the
exact linearization by state coordinate transformation combined with
feedback (further abbreviated as SFB linearization). Our main goal is
to identify the main difficulties of this approach and to consider
applicability of the methods known from the deterministic systems.

The task of SFB linearization is following: given a dynamical
systems~$\Theta$ we are looking for a combination of coordinate
transformation~$\sctt$ and feedback~$\feedbackab$ which will make the
resulting system~$\feedbackab \compose \sctt (\Theta)$ linear and
controllable. One can also define the feedback-less linearization by
coordinate transformation only (here abbreviated as SCT) or several
variants of the input-output linearization. These variants are not
considered here.
\begin{rem}
The notation for composition of mappings sometimes differs;
right-to-left convention is used here:~$f \compose g(x) \define
f(g(x))$.
\end{rem}

The subject of exact linearization of stochastic controlled dynamical
systems lies on the intersection of three branches of science:
differential geometry, control theory, and the theory of stochastic
processes. Each of them is very broad and it is virtually impossible
to cover all details of their combination. Hence it is necessary to
choose a minimalistic simplified model for our problem and to refrain
from most technical details. {\em We decided to represent dynamical
systems under investigation by triplets of smooth vectorfields and to
concentrate on transformation rules for these triplets\/}. The detailed
interpretation of the vectorfield systems (\ie solvability of
underlying differential equations, properties of flows and
trajectories) will be considered only on an informal, motivational level.

For simplicity, we shall confine all the definitions of geometrical
object to the Euclidean space; we will work in a fixed coordinate
system using explicit local coordinates, which may be considered to be
local coordinates of some manifold. This is mainly because we are
unable to capture all consequences  of the modern,
coordinate-free differential geometry to the stochastic calculus (see
\eg~\citet{kendall86},~\citet{malliavin},~\citet{emery89}). We believe
that this approach is quite satisfactory for the majority of practical
applications.

\subsection{Dynamical systems}

\begin{defn}
In this paper, a stochastic dynamical system~$\Theta \define
 \ssysfgx$ is defined to be a triplet of  smooth and bounded
 vectorfields~$f$,~$g$, and~$\sigma$ defined on an open
 neighborhood~$U$ of a point~$\origi{x} \in \RR^n$.  We usually
 call~$U\in \RR^n$ the \defbf{state space},~$f$ the \defbf{drift
 vectorfield},~$g$ the \defbf{control vectorfield}, and ~$\sigma$ the
 \defbf{dispersion vectorfield}.
\end{defn}

From now on, let's assume that all functions, vectorfields, forms, and
distributions are smooth and bounded on~$U$.

In this paper, we will study almost only SISO systems, but in the case
of stochastic MIMO systems with~$m$ control inputs and~$k$-dimensional noise the
symbols~$g$ and~$\sigma$ stand for~$n\times m$ ($n \times k$
respectively) matrix of smooth vectorfields having its rank equal
to~$m$ ($k$ respectively). The class of all
deterministic~$n$-dimensional dynamical systems with~$m$ inputs will
be called~$\clmisysdet$ and the class of stochastic systems
with~$k$-dimensional noise will be denoted with~$\clmisys$.

Similarly, autonomous deterministic dynamical system
corresponds to a single vectorfield and a controlled
deterministic dynamical system corresponds to a vectorfield pair. It
is obvious that this approach is limited to time invariant, affine
systems.

\begin{rem}
  The acronyms SISO and MIMO are used in the usual meaning even for
  systems without outputs, where the wording ``scalar-input''and
  ``vector-input'' will be appropriate. Stochastic systems with~$m =1$
  and~$k=1$ will be considered SISO.
\end{rem}

The definition may be interpreted as follows: there is a stochastic
 process~$x_t$ defined on~$\RR^n$ which is a strong solution of the
 stochastic differential equation $dx_t = f(x_t)\,dt + g(x_t) u(t)\,dt
 + \sigma(x_t)\,dw_t$, with initial condition~$\origi{x}$,
 where~$u(t)$ is a smooth function with bounded derivatives and~$w_t$
 is a one-dimensional Brownian motion. The differential $dw_t$ is just
 a notational shortcut for the stochastic integral.
 
 Details of the theory of stochastic processes are beyond the scope of
 this article. The reader is referred to \citet{wong84},
 \citet{oksendal}, \citet{sagemelsa}, \citet{malliavin},
 \citet{kendall86}, \citet{karatzas}; the text of \citet{kohlman} is
 freely available on the Internet.
 
 Theory of stochastic processes offers several alternative definitions
 of the stochastic integral, among them the \Ito integral and the
 Stratonovich integral; each of them is used to model different
 physical problems. Consequently there are two classes of differential
 equations and two alternative definitions of a stochastic dynamical
 system --- \Ito dynamical systems defined by \Ito integrals and
 Stratonovich systems defined by Stratonovich integrals.

 Serious differences between these integrals exists but from out point
 of view there is a single important one: {\em the rules for
 coordinate transformations of dynamical systems defined by \Ito
 stochastic integral are quite different from the transformation rules
 which are valid for Stratonovich systems\/}.
 
 The definition of the \Ito dynamical system used by us is formally
 equivalent to the definition of the Stratonovich system; the only
 difference will be in the corresponding coordinate transformation.

If necessary, \Ito and Stratonovich dynamical
 systems will be distinguished by a subscript: $\Theta_I \in
 \clmisysito$ and~$\Theta_S \in \clmisysstrat$.

\begin{rem}
In this paper we will use the adjectives {\em \Ito\/} and {\em
  Stratonovich} rather freely. For example we will speak of
  'Stratonovich linearization' instead of `exact linearization of
  stochastic dynamical system defined by Stratonovich integral'.
\end{rem}

\subsection{Transformations}

Furthermore, we will study two transformations of dynamical systems:
the coordinate transformation~$\sctt$ and the
feedback~$\feedbackab$. The definition of these transformation should
be in accord with their common interpretation. This can be illustrated
on the definition of the \defbf{coordinate transformation of a
deterministic dynamical system} $\sctt: \clmisysdet \to \clmisysdet$
which is induced by a diffeomorphism~$T \colon U \to \RR^n$ between two
coordinate systems on an open set~$U \subset \RR^n$.
The mapping~$\sctt$  is defined by:
\begin{align}
 \label{eq:83}
\sctt \dsys{f}{g}{x} \define  \left(
\tantra{T} f , \tantra{T} g ,T(U),T(\origi{x}) 
\right) 
.\end{align}

Recall that the symbol~$\tantra{T}$ stands for the contravariant
transformation~$(\tantra{T} f)_i = \ssum{j}{0}{n} f_j
\parby{T_i}{x_j}$. Moreover, we will require
that the coordinate transformation~$T$ preserves the equilibrium state
of the system \ie~$T(\origi{x}) = 0$.

The definition captures the contravariant transformation
rules for differential equations known from the basic calculus.

Note that the words ``coordinate transformation'' are used in two
different meanings: first as the
diffeomorphism~$T\colon U\to \RR^n$ between coordinates;
second as the mapping between systems~$\sctt: \clmisysdet \to
\clmisysdet$.

Coordinate transformation of stochastic systems distinguish between
\Ito and Stratonovich systems.  One of the major complications of the
linearization problems for \Ito systems is the second-order term in
the transformation rules for \Ito systems:
\begin{defn} \label{def:ctito}
 Let~$U\in \RR^n$ be an open set and let~$T\colon U\to \RR^n$ be a
 diffeomorphism from~$U$ to~$\RR^n$ with bounded
 first derivative on $U$ such that~$T(\origi{x})=0$. The mapping
 $\sctt\colon \clmisysito \to \clmisysito$ will be called a
 \defbf{coordinate transformation of an \Ito dynamical system} induced by
 diffeomorphism~$T$ if the systems $\Theta_1
 \define \ssys{f}{g}{\sigma}{x}$ and $\Theta_2 \define
 \left(\tilde f,\tilde g,\tilde \sigma,T(U),\origi{x}\right)$; $\Theta_2
 = \sctt\left(\Theta_1\right)$ are related by:
 \begin{align}
 \label{eq:86}
 \tilde f &= \tantra{T}f + \ito{\sigma}{T}  \\
 \label{eq:87}
 \tilde g_i &= \tantra{T}g_i \firstfor 1 \le i \le m\\
 \label{eq:89}
 \tilde \sigma_i &= \tantra{T}\sigma_i  \firstfor 1 \le i \le
 k
.\end{align}
\end{defn}

The symbol~$\ito{\sigma}{T}$ stands for the \defbf{\Ito term} which is
a second order linear operator defined by the following relation for
the~$m$-th component of~$\ito{\sigma}{T}$, $1\le m\le n$
\begin{equation}
\label{eq:16}
\ito{\sigma} T_m \define \frahalf \ssum{{i,j}}{1}{n}
\frac{\partial^2 T_m}{\partial x_i x_j} \ssum{l}{1}{k} \sigma_{il}
\sigma_{jl}
.\end{equation}

The transformation rules for Stratonovich
systems~$\sctt\colon \clmisysstrat \to \clmisysstrat$,
$(f,g,\sigma,U,\origi{x}) \mapsto (\tantra{T}f, \tantra{T}g, \tantra{T}\sigma,
T(U), T(\origi{x}))$ are 
equivalent to rules valid for the deterministic systems; only the
rule\eqref{eq:89} for the drift vectorfield must be added.

The difference between the coordinate transformation of \Ito and
Stratonovich systems should be emphasized: in the Stratonovich case
all the vectorfields transform contravariantly; on the other hand, in
the \Ito case, the \Ito term~$\ito{\sigma}{T}$ is added to the drift
vectorfield  of the resulting system.

\psfig{fig:introsfee}{Regular State Feedback}{introsfee} Another
important transformation of dynamical systems is the regular feedback
transformation. A feedback transformation is determined by two smooth
nonlinear functions $\alpha \colon \RR^n \to \RR^m$
and~$\beta \colon \RR^n \to \RR^m \times \RR^m$ defined on~$U$
with~$\beta $ nonsingular for every~$x \in U$ (see
Figure~\ref{fig:introsfee}). Usually, $\alpha $ is written as a column~$m \times
1$ matrix; $\beta $ as a square $m \times m$ matrix.

\begin{defn} \label{def:feedback}
  Let~\dsdefmiuu be a stochastic dynamical system. A \defbf{regular
  state feedback} is the transformation $\feedbackab \colon \clmisys
  \to \clmisys~$, $ (f,g,\sigma,U,\origi{x}) \mapsto \ssys{f + g\alpha
  }{g\beta }{\sigma}{x} $.
\end{defn}

A new input variable~$v$ is introduced by the
relation~$u=\alpha +\beta v$.  Given the feedback~$\feedbackab$
with nonsingular~$\beta $, we can always construct an inverse
relation~$\feedback{a}{b} \define \feedbackab^{-1}$ such
that~$\feedbackab \compose \feedback{a}{b} = \feedback{a}{b} \compose
\feedbackab$ is the identity. The coefficients are related as
follows:~$\beta  = b^{-1}$, $\alpha  = - b^{-1}a$, and~$a = -\beta ^{-1}
\alpha $.

This definition of feedback transformation can be used also for
deterministic systems provided that the drift vectorfield~$\sigma$ is
assumed to be zero.

The symbol~$\combinedtab$ is used to indicate the combination of
 coordinate transformation with feedback~$\combinedtab \define 
 \feedbackab\compose\sctt$ can be interchanged.

\begin{rem}
\label{rem:orderinv}
Observe that the order of feedback and
coordinate transformation in the
composed transformation~$\combinedtab \define \feedbackab \compose
\sctt$
\begin{multline}
 \combinedtab = \sctt \compose \feedbackab \bigdsys{f}{g}{x} =
 \sctt \bigdsys{f+g\alpha }{g\beta }{x}\\
 = \bigdsys{\tantra{T} f + \tantra{T} g\alpha }{\tantra{T} g\beta }{z}
 = \bigdsys{\tantra{T} f + (\tantra{T} g) \alpha  }{ (\tantra{T} g)
 \beta  }{z} \\
= \feedback{\alpha '}{\beta '} \compose \sctt
.\end{multline}
The functions~$\alpha(z)'$, $\beta(z)'$ are equal to~$\alpha(x) $ 
and~$\beta(x) $ written in the $z$~ coordinates $\alpha(z) '=\alpha(x) 
\compose T^{-1}(z)$, $\beta(z) '=\beta(x) \compose T^{-1}(z)$. 
\end{rem}

%======== Subsection ================================================
\subsection{Linearity}

\label{sub:linearity}

The definition of linearity is straightforward in the deterministic
case. In contrast, the stochastic case is more complex, because there
are two ``input'' vectorfields and thereby several degrees of
linearity can be specified.

\begin{defn} \label{def:linear} The deterministic dynamical
  system~$\Theta_D = (f,g,U,0) \in \clmisysdet$ is \defbf{linear} if
  the vectorfield~$f$ is a linear mapping without no constant term and
  the vectorfields~$g_i$ are constant; that is they can be written as
  $f(x)= Ax$, $g(x)= B$ with~$A$ a square~$n \times n$ matrix and~$B$
  an~$n \times m$ matrix. The matrices must be constant on whole~$U$.
\end{defn}
\begin{defn} \label{def:linearsto}
  The stochastic dynamical system~$\Theta=(f,g,\sigma,U,0)$ is:
\begin{itemize}
\item \defbf{$g$-linear} if the mapping~$f(x) = Ax$ is linear without 
  constant term and~$g(x) = B$ is constant on~$U$.
\item \defbf{$\sigma$-linear} if the mapping~$f(x) = Ax$ is linear
  without constant term and~$\sigma(x) = S$ is constant on~$U$.
\item \defbf{$g\sigma$-linear} if it is both~$g$-linear and
  $\sigma$-linear.
\end{itemize}
The matrices~$A$ and~$B$ have the same dimensions as in
Definition~\ref{def:linear}; $S$ is an~$n \times k$ matrix.
\end{defn}

\begin{rem}
 The vectorfield~$g$ is \defbf{constant} on~$U$ if the value
 of~$g(x)$ is the same for every~$x \in U$. The vectorfield~$f$
 on~$U$ is \defbf{linear without constant term} if~$f(\origi{x})=0$
 and the superposition principle~$f(x_1+x_2)=f(x_1)+f(x_2)$
 holds for every~$x_1$, $x_2 \in U$.
\end{rem}

We study systems at equilibrium \ie we require that~$f(\origi{x})=0$
and that all transformations preserve the equilibrium:~$T(\origi{x})
=0$,~$\alpha (\origi{x})=0$, and~$\beta(\origi{x})$ is nonsingular.
The \Ito systems require an additional
condition~$f(\origi{x})+\corr{\sigma}{\origi{x}}=0$.  The non-equilibrium case can
be easily handled by extending the linear model with a constant term.

Moreover we require that the resulting linear systems are
controllable.  A controlled deterministic dynamical linear
system~$\Theta_D = (Ax,B,\RR^n,0) \in \clmisysdet$ is
\defbf{controllable} if its first $n$ repeated brackets form an
$n$-dimensional space
\begin{equation}
 \label{eq:68}
 \dim \left\{A^kB, 0\le k\le n-1 \right\} = n
.\end{equation}
Other definitions of controllability of linear
systems exist. For example Theorem 3.1 of~\citet{zhou} gives six
definitions with proofs of equivalence.

The controllability property deserves some attention in the stochastic
case. The linear stochastic dynamical system is characterized by two input
vectorfields $g(x)=B$ and $\sigma(x)=S$. 
\begin{enumerate}
\item The definition of the controllability for the control
  vectorfield~$g$ is identical to the deterministic case;
  \ie~\eqref{eq:68} must be satisfied. This property will be called
  \defbf{$g$-controllability}.
\item We will also define \defbf{$\sigma$-controllability} as the
  requirement that the repeated brackets $S,AS,A^2S,\dots,A^{n-1}S$
  form an $n$-dimensional space.
\item Finally, the linear system is \defbf{$g\sigma$-controllable} if
\begin{equation}
\label{eq:3}
 \dim \left( \left\{A^kB, 0\le k\le n-1 \right\} \bigcup \left\{A^kS, 0\le
     k\le n-1 \right\} \right)= n
.\end{equation}
\end{enumerate}
In this paper, we do not deal with the reachability, controllability,
accessibility, observability, and similar properties of nonlinear
systems.

\begin{defn} \label{def:itogssisosfb}
  Let~\dsdefmiuu be a dynamical system such that~$f(\origi{x}) =0$. We
  call the combination of a coordinate transformation~$\sctt$ and a
  regular feedback~$\feedbackab$ such that $T(\origi{x})=0$,
  $\alpha (\origi{x})=0$, and $\beta (\origi{x})$ is nonsingular the
  \defbf{linearizing transformation } of $\Theta$ if the
  transformation~$\feedbackab \compose \sctt $ converts~$\Theta$ into
  a~{\em controllable\/} linear system.
 
For stochastic system we distinguish:
\begin{romanenu}
\item \defbf{$g$-linearizing transformation} which transforms~$\Theta$ into
  a~$g$-linear and~$g$-controllable system
\item \defbf{$\sigma$-linearizing transformation} which transforms~$\Theta$ into
  a~$\sigma$-linear and~$\sigma$-controllable system 
\item \defbf{$g\sigma$-linearizing 
 transformation}  which transforms~$\Theta$ into
  a~$g\sigma$-linear and~$g$-controllable system 

Note that for $g\sigma$-linearization we
require~$g$-controllability. This is slightly stricter requirement
than~$g\sigma$ controllability but it should be naturally fulfilled by
the majority of practical control systems. This requirement cancels
many ``uncomfortable'' linear forms. Consider for example the system
with prefilter of Figure~\ref{fig:prefilter} which
is~$g\sigma$-controllable but~$g$-uncontrollable.
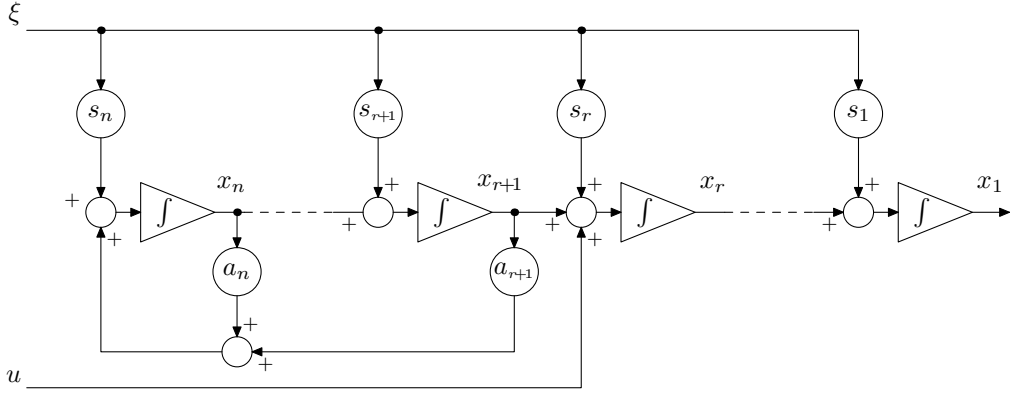
\end{romanenu}\psfig{fig:prefilter}{Dynamical System with
a Prefilter.}{prefilter}
 The system~$\Theta$ is \defbf{linearizable} if there exists 
linearizing transformation of~$\Theta$.
\end{defn}

\subsection{Computational Issues}

In most practical circumstances, computational issues are the limiting
factor of any application of differential geometric methods in
control.

The equations of exact linearization algorithms must be dealt in a
symbolic form. Even the simplest exact linearization problems are
extremely complex from the computational point of view.  Therefore,
the computer algebra tools are often employed. The results presented
in this paper were tested by the author on  few simulations of
control systems in the symbolic  system \Mathematica.

Of course, the computer algebra has apparently serious limitations and
drawbacks. Viability of the symbolic computational approach to the
problems of the nonlinear control is studied by \citet{jager95}. Some
very useful theoretical notes about the symbolic computation can be
found in~\citet{winkler}. Unfortunately, the limited scope
of this article does not allow deeper discussion of these subjects.

\subsection{Applications}

We propose, very briefly, two applications of the theory presented here: 
\begin{romanenu}
\item \defbf{Control ---}
 a dynamical systems~$\Theta$ obtained by exact
linearization will be controlled using the linear
feedback law:
\begin{equation}
 \label{eq:141}
 v = Kz + \kappa \nu
,\end{equation}
where~$K$ is a row matrix of feedback gains,~$\kappa$ is an input
gain, $z$~is the state vector, $v$~is the original control input,
 and~$\nu$ is a new control input. 

Two approaches can be studied --- classical linear control methods and
the more sophisticated stochastic optimal control  approach studied
for example by \citet{oksendal}.

The~$g\sigma$-linear systems are natural candidates for such
approach because  the other linear forms leave certain
part of the resulting system nonlinear.

\item \defbf{Filtering} --- the filtering problem is probably the most useful
application of the theory of stochastic processes. We want
to give the best estimate to the state of a dynamical system defined
by the stochastic differential equation:
\begin{align}
\label{eq:118} dx_t = f(x_t)\,dt + \sigma_f(x_t)\,dw_tf;
\end{align} based on observations of the from:
\begin{align}
\label{eq:119}y_t = h(x_t) +
\sigma_h(x_t)\,dw_th .\end{align} $x_t$ is an $n$-dimensional stochastic process, $f$, $\sigma_f$,
$\sigma_h$ are smooth vectorfields and $h$ is a
smooth function.

It would be interesting to use exact linearization of the nonlinear
system to design an exact linear filter. Unfortunately, this idea has
no direct association with the linearization results presented below,
because it requires {\em output\/} exact linearization or
linearization of an autonomous system. Therefore, it would be helpful
to extend our results to these cases in the future.

\end{romanenu}

\subsection{Previous Work}
\label{ssub:laal}

The problem of SFB~$g$-linearization of SISO dynamical system defined in the
 \Ito formalism has been 
studied by~\citet{lahdhiri}.  The authors derive
equations corresponding to~\eqref{eq:50},~\eqref{eq:51} (eq 14, 15, 16
in~\citez{lahdhiri}). These equations are combined and then
reduced to a set of PDEs of a single
unknown function~$T_1$. Because there is no commuting relation
similar to~\eqref{eq:103} the equations contain partial derivatives
of~$T_1$ up to the~$2n$-th order (eq 23 in~\citez{lahdhiri}).
Next, the authors propose a lemma (Lemma 1) that identifies the
linearity conditions with non-singularity and involutivness
of~$\distrofg{n-2}$. Unfortunately, we disagree with this result.

 It can be easily verified that for~$\sigma=0$ this
statement does not correspond to the deterministic conditions
(Proposition~\ref{prop:d1sfb2}), because the deterministic case
requires non-singularity up to the~$(n-1)$-th bracket, not only up to
the~$(n-2)$-th one. Second, although the method of finding~$T_1$
was given (solving PDE), we do not think
that the existence of~$T_1$ was proved.

After this paper was finished, we discovered recent works
of~\citet{pan02} and~\citez{pan01}. In the article~\citez{pan01}
Pan defines and solves the problem of {\em feedback complete
linearization of stochastic nonlinear systems}. In our terminology,
this problem is equivalent to SFB MIMO input--output \Ito $g\sigma$
linearization which was not studied by us. 

In~\citez{pan02} Pan declares and proves so called \index{invariance
  under transformation rule}{\em invariance under transformation rule}
  which is exactly equivalent to our Theorem~\ref{prop:corr} which is
  probably the most important result of our paper.

Althougth the problems solved by Pan were slightly different, he uses
the same equivalence --- Theorem~\ref{prop:corr}.  This proves that
our conclusions about applicability of the \index{correcting
term}correcting term are perfectly valid.

In~\citez{pan02} Pan
consider three other \index{canonical form}canonical forms of
stochastic nonlinear systems, namely the \index{noise-prone strict
feedback form}noise-prone strict feedback form, \index{zero dynamics
canonical form}zero dynamics canonical form and \index{and observer
canonical form}observer canonical form also not studied by us. 

%%%%%%%%%%%%%%% Section %%%%%%%%%%%%%%%%%%%%%%%%%%%%%%%%%%%%%%%%%%%%%
\section{Deterministic Case}
\label{sec:detcase}

In this section we recapitulate the results of the SFB and SCT exact linearization
theory for SISO systems. For detailed treatment and proofs we refer to
existing literature, above all  the
classical monographs of~\citet{isidori85} and~\citet{nijmeijer94}. For a
very readable introduction to the field we refer to the seventh
chapter of~\citet{vidyasagar93}. The books also contain extensive
bibliography. The monograph of~\citet{isidori85} builds mainly on the
concept of relative degree. In contrast we will emphasize the approach
of~\citet{vidyasagar93} because the method is  more suitable for the
stochastic case.

\subsection{Useful Relations}

The solution of the SFB linearization problem as presented
here  uses  the Leibniz rule
\begin{align}
\label{eq:63}
\lie{\lbrack f,g\rbrack }{\alpha } =
\lie{f}{(\lie{g}{\alpha })} - \lie{g}{ (\lie{f}{\alpha })}
\end{align} 
with~$f,g$ smooth vectorfields on~$U$; $\alpha\colon U\to \RR $ is
a smooth function. The recursive form of the Leibniz rule allows to
simplify the chains of differential equations for the
transformation~$T$. This can be expressed in the form of the following
statement:

For all $x \in U$, $k \ge 0$ these two sets of conditions are
equivalent:

\begin{align}
\label{eq:2}
\text{(i)}&\qquad&\lie{g}{\alpha } = \lie{g}{\lie{f}{\alpha }} = 
\cdots = \lie{g}{\multilie{f}{k}{\alpha }}=0 \\
\label{eq:10}
\text{(ii)}&&\lie{g}{\alpha } = \lie{\adfg{}}{\alpha } = \cdots =
\lie{\adfg{k}}{\alpha }=0 .\end{align}

Recall that the symbol~$\lie{f}{h}$ stands for the Lie derivative 
 defined by~$\lie{f}{h} = \langle f,\grad\,h \rangle = \ssumi{n} f_i(x)
 \parby{}{x_i} h(x)$. Higher order Lie derivatives
 can be defined recursively
$\multilie{f}{0}{h} = h$, $\multilie{f}{k+1}{h} =
 \lie{f}{\multilie{f}{k}{h}}$ for~$k \ge 0$. The Lie Bracket is defined as
~$[f,g] \define \parbyx{g}f-\parbyx{f}g$; there is also a recursive
 definition:
\begin{equation}
 \ad{f}{0}{g} \define g; \quad \ad{f}{k+1}{g} \define
 \left[f,\ad{f}{k}{g}\right] \quad \text{for } k \ge 0 
.\end{equation}

Another very important result of the differential geometry is
invariance of the Lie bracket under the tangent
transformation~$\tantra{T}$ 
(see~\citet{nijmeijer94} Proposition 2.30 \pg~50):

Let~$T\colon U\to \RR^n$ be a diffeomorphic coordinate transformation, and~$f$ and~$g$
be smooth vectorfields. Then
\begin{align}
 \label{eq:103}
 \tantra{T} [f,g] = [\tantra{T} f, \tantra{T} g] .\end{align}

\subsection{SFB Linearization}

Every controllable linear system may be, by a linear coordinate
transformation, transformed to the controllable canonical form
(\citet{kalman}). Furthermore, this controllable canonical form can
be always transformed into the integrator chain by a linear
regular feedback.  Therefore, the integrator chain is a canonical
form for all feedback linearizable systems. (See \citet{vidyasagar93}
section 7.4). Consequently, the  equations of the integrator chain can be compared
with the equations of the nonlinear systems and the following
proposition can be proved:

\begin{prop} \label{prop:d1sfb1a}
 There is a SFB linearizing transformation~$\combinedtab$ of a SISO
 deterministic dynamical system~\dsdef into a controllable linear
 system if and only if there is a solution~$T_1, T_2, \dots,
 T_n\RntoR$ to the set of partial differential equations
 defined on~$U$
\begin{alignat}{2}
\label{eq:150}
\lie{f}{}{T_i} &= T_{i+1} \firstfor 1 \le i \le {n-1}\\
\label{eq:151}
\lie{g}{}{T_i} &= 0 \nextfor 1 \le i \le {n-1}\\
\label{eq:152}
\lie{g}{}{T_n} &\ne 0 .\end{alignat} Then the
feedback is defined as follows:
\begin{align}
\label{eq:57}
\alpha &=-\frac{\lie{f}{T_n}}{\lie{g}{{T_n}}}\qquad\qquad
\beta =\frac{1}{\lie{g}{T_n}} .\end{align}
\end{prop}
\begin{proof}
See \citet{vidyasagar93} equations 7.4.20--21.
\end{proof}

\begin{prop} \label{prop:d1sfb1}
 The SFB linearizing transformation~$\combinedtab$ of a SISO
 deterministic dynamical system~\dsdef into a controllable linear
 system  exists if and only if there is
 a solution~$\lambda \RntoR$ to the set of partial differential
 equations:
\begin{alignat}3
\label{eq:8}
\biglie{\adfg{i}}{\lambda }&=0\firstfor{\iva{n-2}}\\
\label{eq:9}
\biglie{\adfg{{n-1}}}{\lambda }&\not=0 .\end{alignat} The
linearizing transformation~$T(x)$ is given by:
\begin{alignat}3
\label{eq:11}
T_{i}&=\multilie{f}{{i-1}}\lambda  \firstfor 1\le i\le n\\
\alpha &=\frac{-\multilie{f}{{n}}\lambda 
}{\lie{g}{}\multilie{f}{{n-1}}\lambda }&\qquad\qquad&
\label{eq:12}
\beta =\frac{1}{\lie{g}{}\multilie{f}{{n-1}}\lambda }
.\end{alignat}
\begin{proof}
See \citet{vidyasagar93} equations 7.4.23--33 and \citet{nijmeijer94} Corollary 6.16.
\end{proof}

\end{prop}
Finally, the  geometrical conditions for the
existence of the linearizing transformation  are studied.

\begin{theo} \label{prop:d1sfb2}
 A deterministic SFB linearizing transformation of~\dsdef into a
 controllable linear system exists if and only if the
 distribution~$\Delta_{n} \define \sspan\left\{{\adfg{i},
 \iva{n-1}}\right\}$ is nonsingular on~$U$ and the
 distribution~$\Delta_{n-1} \define \sspan\left\{{\adfg{i},
 \iva{n-2}}\right\}$ is involutive on~$U$.
\end{theo}
\begin{proof}
See \citet{nijmeijer94} Corollary 6.17, \citet{vidyasagar93} Theorem
7.4.16, \citet{isidori89} Theorem 4.2.3 .
\end{proof}

\subsection{SCT Linearization}

\begin{theo} \label{prop:s1ctt2t}
  There is a SCT~linearizing transformation~$\sctt$ of a deterministic
  MIMO system~\dsdefmi into a controllable linear system if  and only
  if there exists a reordering of the vectorfields~$g_1 \dots g_m$ and
  an~$m$-tuple of integers~$\kappa_1 \le \kappa_2 \le \dots \kappa_m$
  with~$\ssumi{m} \kappa_i = n$ called the \defbf{controllability
    indexes} such that the following conditions are satisfied for
  all~$x \in U$:
 \begin{align}
 \label{eq:43}
 \text{(i)}&\qquad\dim\left(\sspan\left\{( \ad{f}{j}{g_i}(x),
 \iva{m}, \jva{\kappa_i-1})\right\}\right) = n\\
 \label{eq:44}
 \text{(ii)}&\qquad[\ad{f}{k}{g_i},\ad{f}{l}{g_j}]=0
 \qquad\text{for}\qquad 0 \le k+l \le \kappa_i+\kappa_j-1,\, 0 \le
 i,j \le m . \end{align}
\end{theo}
\begin{proof}
See  \citet{nijmeijer94} Theorem 5.3 and Corollary 5.6.
\end{proof}

The following corollary can be verified for SISO systems:

\begin{coroll} \label{prop:s1ctt2a}
 For a SISO system with~$m=1$ the
 condition (ii) of~\ref{prop:s1ctt2t} can be simplified as
 follows:
\begin{align}
\label{eq:45}
[g,ad^l_f g] = 0, \quad l = 1, 3, 5, \dots, 2n-1, \forall x \in U
.\end{align}
\end{coroll}
\begin{proof}
See  \citet{nijmeijer94} Corollary 5.6 and the text which follows.
\end{proof}

%%%%%%%%%%%%%%% Section %%%%%%%%%%%%%%%%%%%%%%%%%%%%%%%%%%%%%%%%%%%%%
\section{Transformations of \Ito Dynamical Systems}
\label{sec:stotrans}

The transformation rules of \Ito systems are motivated by the  \Ito
differential rule (see \eg~\citet{wong84} Section 3.3), which defines
the influence of nonlinear coordinate transformations on  \Ito stochastic
processes.

The \Ito differential rule applies to the situation where a scalar
valued stochastic process~$x_t$ defined by a stochastic differential
equation~$dx_t = f(x_t)\,dt + \sigma(x_t)\,dw_t$ with~$f\RtoR$ and~$\sigma\RtoR$
smooth real functions and~$w_t$  a Brownian motion, is transformed
by a a diffeomorphic coordinate transformation~$T\colon \RR \to \RR$.
Then the stochastic process~$z_t\define T(x_t)$ exists and is an \Ito
process. Further, the process~$z_t$ is the solution of the stochastic
differential equation
\begin{equation} \label{eq:126}
dz_t = \parbyx{T} {f(x_t)}\,dt + \parbyx{T} {\sigma(x_t)}\,dw_t + \itoxs\,dt
.\end{equation}
All details together with a proof are
available for example in~\citet{karatzas}.

The \Ito rule can be also derived for the multidimensional case: for
the~$m$-th component of an~$n$-dimensional stochastic process the \Ito
rule can be expressed as follows:
\begin{multline}
 dz_m = \ssumi{n} \parby{T_m}{x_i} f_i\,dt + \frahalf \ssum{i}{1}{n}
 \ssum{j}{1}{k} \parby{T_m}{x_i} \sigma_{ij}\,dw_j + \frahalf
 \ssum{{i,j}}{1}{n} \frac{\partial^2 T_m}{\partial x_i x_j}
 \ssum{l}{1}{k} \sigma_{il}
 \sigma_{jl}\,dt \\
 =\lie{f}T_m\,dt + \ssum{j}{1}{k} \lie{\sigma_j}{T_m}\,dw_j +
 \ito{\sigma}T_m\,dt
.\end{multline}

For the most common case with scalar noise~$k=1$ the equation
 can be further simplified to:
\begin{equation}
dz_m =\lie{f}{T_m}\,dt + \lie{\sigma}{T_m}\,dw +
\frahalf \ssum{{i,j}}{1}{n}
\parbyxx{T_m}{x_i}{x_j} \sigma_{i} \sigma_{j}\,dt 
.\end{equation}
The operator~$\ito{\sigma} T_m$ is sometimes written using
matrix notation as:
\begin{equation}
\ito{\sigma} T_m = \frahalf \trace \left(\sigma^T\sigma \parbysec{T_m}{x}\right)
.\end{equation}
Generally, $\ito{\sigma}$ vanishes for linear~$T$ or zero~$\sigma$.

%======== Subsection ================================================
\subsection{The Correcting Term}
\label{sub:corr}

In this section we introduce an extremely useful equivalence
between \Ito and Stratonovich systems, which allows to use some
Stratonovich linearization techniques for \Ito problems. The
motivation is following: let~\isdefmi be an  \Ito system. We
are looking for a Stratonovich system $\Theta_S = \ssys{\alterbar
 f}{\alterbar g}{\alterbar \sigma}{x}$ such that the trajectories
of~$\Theta_I$ and~$\Theta_S$ are identical. The aim is to find
equations relating the quantities~$\alterbar f$,~$\alterbar
g$, and~$\alterbar \sigma$ with~$f$,~$g$, and~$\sigma$.

\begin{defn} \label{def:corr}
 Let~$\Theta_{1I} = \ssys{f}{g}{\sigma}{x} \in \clmisysito$ be
 an~$n$-dimensional \Ito dynamical
 system with~$k$-dimensional Brownian
 motion~$w$. The vectorfield $\corr{\sigma}{x}$ whose ~$r$-th
 coordinate is equal to
\begin{alignat}3 
\label{eq:5}
 (\corr{\sigma}{x})_r &= -\frahalf \ssum{i}{1}{n} \ssum{j}{1}{k}
 \parby{{\sigma_{rj}}}{{x_i}} \sigma_{ij} \firstfor 1 \le r \le n
\end{alignat}
is called the \defbf{correcting term}.
Note that the derivative is always evaluated in the corresponding
coordinate system.
Further, let us define the \defbf{correcting mapping} $\Corr{\sigma} \colon \clmisysito \to \clmisysstrat$ by
\begin{align}
\label{eq:21}
 \Corr{\sigma} (f,g,\sigma,U,\origi{x}) \define
 (f+\corr{\sigma}{x},g,\sigma,U,\origi{x}) .\end{align}
 \end{defn}
 The general treatment of the subject can be found for example
in~\citet{wong84} \pg~160 or in \citet{sagemelsa}.  The following
theorem describes the
behavior of the correcting term under the coordinate transformation.
\begin{theo} \label{prop:corr}
  Let~\isdef be a one-dimensional \Ito dynamical system. Let $T$ be a
  diffeomorphism defined on~$U$ and the symbols~$\scti{T}$ and
  $\scts{T}$ denote a \Ito coordinate transformation and a
  Stratonovich coordinate transformation induced by the same
  diffeomorphism~$T$ and~$\tilde\sigma = \tantra{T}{\sigma}$.  Then the following diagram commutes:
\begin{equation}
\begin{CD}
\label{eq:4}
\Theta_{1I} @>\scti{T}>> \Theta_{2I}\\
@V{\Corr{\sigma}}VV @AA{\Corr{{\tantra{T} \sigma}}^{-1}}A\\
\Theta_{1S} @>\scts{T}>> \Theta_{2S}\\
\end{CD}
.\end{equation}
In other words:
 \begin{align}
 \label{eq:19}
  \scti{T} &= {\Corr{{\tilde\sigma}}{}}^{-1} \compose \scts{T} 
  \compose \Corr{\sigma}{} \quad\text{and} \\
 \label{eq:20}
  \scts{T} &=  {\Corr{\sigma}{}}^{-1}  \compose \scti{T}
  \compose \Corr{{\tilde\sigma}}{}  
 . \end{align}

The notation~$\Corr{{\sigma}}{}^{-1}$ is used to denote the
inverse mapping 
\begin{align}
\label{eq:23}
 {\Corr{\sigma}{}}^{-1} (f,g,\sigma,U,\origi{x}) \define
 (f-\corr{\sigma}{x},g,\sigma,U,\origi{x}) 
.\end{align}
\end{theo}

\begin{proof}
 The correcting term~$\corr{\sigma}{x}\RntoRn$
 is equal to
\begin{equation}
\corr{\sigma}{x} = -\corrx
.\end{equation}
The transformations specified in the diagram~\eqref{eq:4} will be
evaluated in the following order:
\begin{equation}
\begin{CD}
\label{eq:6}
\Theta_{1I} @>({\operatorname a})>> \Theta_{2I} \\
@V({\operatorname b})VV \\
\Theta_{1S} @>({\operatorname c})>> \Theta_{2S}
@>({\operatorname d})>> \Theta_{3I}\\
\end{CD}
\end{equation}
We want to prove the equivalence of~$\Theta_{2I}$ and~$\Theta_{3I}$.
The symbols (a) and (c) denote \Ito
coordinate transformations; the
symbols (b) (d) stand for the correcting
mapping and its inverse. Note that the systems~$\Theta_{2I}$,
$\Theta_{2S}$, and $\Theta_{3I}$ are defined in the~$z$-coordinate
systems. Further, let
\begin{align}
 \label{eq:243}
 \tilde\sigma&\define\trav{\sigma}\\
 \kappa&\define \tantra{T} \corr{\sigma}{x} = -\trav{(\corrx)}\\
 \origi{z}&\define T(\origi{x}) .\end{align} Then
\begin{align}
 \text{(a)}&\qquad \Theta_{2I} =\bigssys{\trav{f} + \itoxs }{
 \trav{g}}{ \tilde\sigma }{z}\\
 \text{(b)}&\qquad \Theta_{1S} =\bigssys{f - \corrx }{ g }{
 \sigma }{x} \\
 \text{(c)}&\qquad \Theta_{2S} =\bigssys{\trav{\left(f-\corrx\right)}
 }{ \trav{g} }{
 \trav{\sigma} }{z}\\
 \text{(d)}&\qquad \Theta_{3I} =\bigssys{\trav{f} + \kappa + \frahalf
 \parbyz{\tilde\sigma}{\tilde\sigma} }{ \trav{g} }{ \tilde\sigma
 }{z}
.\end{align}

All the terms in (a) are equivalent to the respective terms in (d)
except for the drift terms containing functions of~$\sigma$. Therefore,
we continue comparing these terms only. For (a):

\begin{align}
 \label{eq:37}
 L &\define \itoxs .\end{align} For (d):
\begin{multline}
\label{eq:38}
R \define \kappa + \corrxx{{\tilde\sigma}}{z} = \kappa + \frahalf
\parbyz{} \left( \parbyx{T} \sigma \right) \parbyx{T}\sigma = \kappa +
\frahalf \parbyz{x} \parbyx{} \left( \parbyx{T} \sigma \right)
\parbyx{T}\sigma = \\
\kappa + \frahalf \left( \parbysec{T}{x}\sigma +
 \parbyx{T}\parbyx{\sigma} \right) \sigma = \kappa + \frahalf
\parbysec{T}{x}\sigma^2 - \kappa= \itoxs .\end{multline} Thus~$L=R$
 and $\Theta_{2I} = \Theta_{3I}$.
\end{proof}

Theorem~\ref{prop:corr} is valid also for combined
 transformations:
\begin{coroll} \label{prop:corrcorol}
  Let~\isdef, $T$, $\scti{T}$, and
  $\scts{T}$ have the same meaning as in Theorem~\ref{prop:corr}. 
 Then the following diagram commutes for arbitrary regular feedback~$\feedbackab$:
\begin{equation}
\begin{CD}
\label{eq:4000}
\Theta_{1I} @>\scti{T}>> \Theta_{2I} @>\feedbackab>> \Theta_{4I}\\
@V{\Corr{\sigma}}VV @AA{\Corr{{\tantra{T} \sigma}}^{-1}}A @AA{\Corr{{\tantra{T} \sigma}}^{-1}}A\\
\Theta_{1S} @>\scts{T}>> \Theta_{2S}  @>\feedbackab>> \Theta_{4S}\\
\end{CD}
.\end{equation}
\end{coroll}

\begin{proof}
  We want to prove equivalence of~$\Theta_{4I}$
  and~${\Corr{{\tantra{T} \sigma}}^{-1}} \Theta_{4S}$. 
  
  The control and dispersion vectorfields of~$\Theta_{4I}$
  and~$\Theta_{4S}$ are identical and they are not influenced by the
  correcting mapping.

  Using the notation of Theorem~\ref{prop:corr} we can express the
  drift term of~$\Theta_{4I}$ as~$\tantra{T}f + L + g\alpha $. The
  drift term of~${\Corr{{\tantra{T} \sigma}}^{-1}} \Theta_{4S}$ 
  is~$\tantra{T}f + R + g\alpha  $. The effect of feedback is purely
  additive and both the systems are equal.
\end{proof}

At first glance the correcting term is rather
surprising. How can the second derivative of~$T$ in (\ref{eq:37}) be
compensated by the correcting term, which does
not contain the~$T$ at all? The answer is quite simple: the second
derivative is hidden in the correcting term
implicitly because the correcting term depends
on the coordinate system in which the system~$\Theta_{1I}$ is defined.
The derivative~$\parbyx{\sigma}$ contained in the
correcting term is always taken in the
appropriate coordinate system. To emphasize the dependence of the
correcting term on the coordinate system, we
will never omit the independent variable (\eg $x$ or $z$) from the
symbol~$\corr{\sigma}{x}$.

 Proposition~\ref{prop:corr} is valid for general multidimensional
 systems~\isdefmi; the proof is purely mechanical and is not presented
 here.

Let us now turn our attention to several special cases of the
correcting mapping.
\begin{coroll} \label{prop:corrvectr}
 Let~$\Theta_I = \ssysfgx \in \clssysito(n,m,1)$ be
 an~$n$-dimensional stochastic dynamical system with an
 one-dimensional Brownian motion~$w$.
 The~$r$-th coordinate of the correcting
 term~$(\corr{\sigma}{x})_r$ is equal to
\begin{equation}
\label{eq:214}
(\corr{\sigma}{x})_r = \frahalf \ssum{i}{1}{n}
\parby{{\sigma_{r}}}{{x_i}} \sigma_{i} = \frahalf
\lie{\sigma}{\sigma_r}
\rfor{1\le r\le n}
.\end{equation}
\end{coroll}
\begin{proof}
 Substitute~$k = 1$ into~\eqref{eq:5}.
\end{proof}
\begin{coroll}
 \label{prop:ic1}
 For systems with one-dimensional noise ($k=1$) define the matrix
 valued \Ito term~$\ito{\sigma}{T}$ for~$T\RntoRn$
 with components~$T_i$, $1\le i\le n$, as a column $n\times 1$
 matrix $\ito{\sigma}{T} \define \left[\ito{\sigma}{T_1},
 \ito{\sigma}{T_2}, \dots,\ito{\sigma}{T_n}\right]^T$.
 
 Then the relations~\eqref{eq:214} can be expressed as
\begin{align}
\label{eq:215}
\ito{\sigma}{T} = \tantra{T} \left(\corr{\sigma}{x}\right) -
\corr{{\tilde\sigma}}{z} .\end{align}
\end{coroll}
\begin{proof}
 The proof is almost identical to that of
 the multidimensional variant of Corollary~\eqref{prop:corr}.The symbols can be
 identified as follows:
 \begin{alignat}2
 L_i &= \left(\ito{\sigma}{T}\right)_i \firstfor 1\le i\le n \\
 \kappa_i &= \left( \tantra{T} \left(\corr{\sigma}{x}\right)\right)_i \\
 R_i &= \left(\tantra{T} \left(\corr{\sigma}{x}\right) \right)_i -
 \left(\corr{{\tilde\sigma}}{z}\right)_i .\end{alignat}
\end{proof}
\begin{coroll} 
 \label{prop:ic2}
 Assume that the conditions of Proposition~\ref{prop:ic1} hold. The
 relation~\eqref{eq:215} can be written as:
 \begin{alignat}3
 \label{eq:140}
 \frahalf \lie{\sigma}{\lie{\sigma}{T_i}} &=
 \ito{\sigma}{T_i} - \lie{\corr{\sigma}{x}}{T_i} \firstfor 1
 \le i \le n
. \end{alignat} 
\end{coroll}
\begin{proof}
 The formula can be expressed as:
 \begin{multline}
 \label{eq:139}
 \frahalf \lie{\sigma}{{\lie{\sigma}{T_i}}} = \frahalf
 \lie{\sigma}{ \left( \ssum{j}{1}{n} \parby{T_i}{x_j} \sigma_j
 \right) } = \frahalf \ssum{l}{1}{n} \sigma_k \parby{}{x_k}
 \left( \ssum{j}{1}{n} \parby{T_i}{x_j} \sigma_j
 \right) =\\
 \frahalf \ssum{k,j}{1}{n} \left( \sigma_j\sigma_k
 \parbyxx{T_i}{x_k}{x_j} +\parby{T_i}{x_j}
 \parby{\sigma_j}{x_k} \sigma_k \right) = \ito{\sigma}{T_i} -
 \lie{{\corr{\sigma}{x}}}{T_i} . \end{multline}
\end{proof}

%======== Subsection ================================================
\subsection{Composition of Coordinate Transformations of \Ito Systems}
\label{sub:stgroup}

The set of all deterministic coordinate transformations~$\sctt$
together with composition~$\sct{RS} \define \sct{S} \compose \sct{R}$
forms a group. Obviously, this fact is a straightforward result of the
behavior of the contravariant transformation and therefore an
analogous statement must hold for Stratonovich systems.  Surprisingly,
this is valid also for \Ito systems as will be shown here. This has an
important consequence: we may always find the inverse transformation
to a given coordinate transformation of \Ito systems. We will prove
the following assertion:
\begin{theo}
 \label{prop:compogr}

 Let~$\scti{R}$, $\scti{S} \in \clsysito$ be
 coordinate transformations of
 one-dimensional \Ito systems  induced by
 diffeomorphisms~$R$ and~$S$. Then
\begin{equation}
 \label{eq:123}
 \scti{S} \compose \scti{R} = \scti{S \compose R} 
.\end{equation}
\end{theo}
\begin{proof}
 We will transform the system in two
 different ways and show that the results are equal.
\begin{enumerate}
\item In the first method the system~$A=(0,0,a,U,\origi{x})$ which
  corresponds to a differential equation~$dx = a(x)\,dw$ will be
  transformed twice:
\begin{enumerate} \item first,  by~$y=R(x)$ to~$y$
 coordinates \item and then the result~$B=(g,0,b,U,\origi{x})$ which
  corresponds to~$dy = g(y)\,dt + b(y)\,dw$
  by~$z=S(y)$ to~$z$ coordinates. 
\end{enumerate}
\item The other method
 transforms the system~$A$ only once by~$z=T(x) = S(R(x)) = (S
 \compose R)(x)$.
\end{enumerate}
Without loss of generality, the equation~$dx = a(x)\,dw$ is assumed
 to have no drift term because the drift term transforms in the
 contravariant fashion. The derivatives will be
 denoted by~$\parbyx{T(x)} \definer T'$,~$\parbyx{R(x)} \definer
 R'$,~$\parby{S(y)}{y} \definer S'$ and similarly for~$T''$,
 $R''$ and $S''$. Note that the prime is always used to denote
 derivatives by the argument of the function.
 
 The transformation by~$R$ gives:
 \begin{align}
 \label{eq:64} 
 dx &= a \, dw\\
 \label{eq:266} 
 dy &= R' a\,dw + \frahalf a^2 R'' \,dt .\end{align} Thus the
 coefficients of the second SDE are defined
 by
\begin{align}
 \label{eq:267}
 b(y) &\define (R'a) \compose R^{-1}(y)\\
 \label{eq:268}
 g(y) &\define (\frahalf a^2 R'') \compose R^{-1}(y) .\end{align}

The second transformation (by~$S$) gives
\begin{align}
\label{eq:269}
dz &= \left( S' g + \frahalf b^2 S'' \right)\,dt + S' b \,dw =\\
\label{eq:270}
&= \left(\frahalf a^2 S' R'' + \frahalf (R')^2 a^2 S'' \right) \,dt +
S' R' a \,dw = \\
\label{eq:272}
&= \frahalf a^2 \left(S'R'' + (R')^2 S'' \right)\,dt + S'R' a\,dw =\\
\label{eq:273}
&= \frahalf a^2 T''\,dt + T' a\,dw .\end{align} The last equality
follows from the fact that~$T''=S''(R')^2+S'R''$.
\end{proof}
One can verify the multidimensional case in the same
spirit. 
%======== Subsection ================================================
\subsection{Invariants}
\label{sub:stoinvar}

In the deterministic case, some useful propositions about the
invariant properties for example the Leibniz rule~\eqref{eq:63} and
the relation~\eqref{eq:2} were employed. 

Unfortunately, we have not found any analogy for the \Ito systems yet.
To point out the main complications, we will analyze the \Ito
equivalent of the Leibniz rule~\eqref{eq:63}, which is essential for
reducing the order of partial differential equations in the
deterministic exact linearization.

If the Lie derivative~$\lie{g}{}$ is interpreted
as a general first order operator
\begin{align}
 \lie{g}{} = \ssumi{n}g_i \parby{}{x_i}
\end{align}
then the commutator of two such first order
operators $\lie{f}{\lie{g}{}} - \lie{g}{\lie{f}{}}$ is also a first
order operator~$\lie{\lbrack f,g\rbrack }{}$ (see~\eqref{eq:63}).

Similarly, define the general second order operator as
\begin{align}
\label{eq:17}
O(g,G) &\define \ssumi{n}g_i \parby{}{x_i} + \ssum{i,j}{0}{n}
G_{ij}\parbyxx{}{x_i}{x_j}
\end{align}
where $g_i,G_{ij}\RntoR$ for $1 \le i,j \le n$. We can compute
the commutator
\begin{align}
\label{eq:213}
C(f,F,g,G) \define O(f,F)O(g,G) - O(g,G)O(f,F)
\end{align} of such
second order operators. If this commutator was also
a second order operator (\ie there were~$\varphi$ and $\Phi$
such that the operator~$C(f,F,g,G) = O(\varphi,\Phi)$), then we would
be able to simplify any PDEs of stochastic
transformations. (See Proposition~\ref{prop:sfbooo}).

Because the operator~$O$ is linear, \ie~$O(f,F) = O(f,0) + O(0,F)$, we
 can split the computation into four independent, reusable parts:
\begin{multline}
 C(f,F,g,G) = O(f,F)O(g,G) - O(g,G)O(f,F) = \\
 =\lie{[f,g]}{} + C(0,F,0,G) + C(f,0,0,G) + C(0,F,g,0)
 .\end{multline}

The first term is already a first order operator. Only the second
 and the third terms need to be computed because the fourth term
 can be obtained from the third one by formal substitution. For
the third term:
\begin{multline}
 O(f,0) O(0,G) = \ssum{i}{1}{n} f_{i} \parby{}{x_i} \left(
 \ssum{k,l}{1}{n} G_{kl}
 \parbyxx{}{x_k}{x_l} \right) = \\
 \ssum{i,k,l}{1}{n} \left( f_{i} \parby{G_{kl}}{x_i}
 \parbyxx{}{x_k}{x_l} + f_{i} G_{kl} \parbyxxx{}{x_i}{x_k}{x_l}
 \right) 
.\end{multline} 
Further,
\begin{multline}
 O(0,G)0(f,0) = \ssum{k,l}{1}{n} G_{kl} \parbyxx{}{x_k}{x_l} \left(
 \ssum{i}{1}{n} f_i \parby{}{x_i}
 \right) =\\
 \ssum{i,k,l}{1}{n} G_{kl} \left( \parbyxx{f_i}{x_k}{x_l}
 \parby{}{x_i} + 2 \parby{f_i}{x_l} \parbyxx{}{x_i}{x_k} + f_i
 \parbyxxx{}{x_i}{x_k}{x_l} \right) .\end{multline}

The intermediate results for the third and the fourth terms
 can be combined into
\begin{multline}
 O(f,0) O(0,G) + O(0,F) O(g,0) - O(g,0) O(0,F) - O(0,G) O(f,0) =\\
 \ssum{i,k,l}{1}{n} \left( f_{i} \parby{G_{kl}}{x_i} - g_{i}
 \parby{F_{kl}}{x_i} + 2 F_{ki} \parby{g_l}{x_i} - 2 G_{ki}
 \parby{f_l}{x_i}
 \right) \parbyxx{}{x_k}{x_l} \\
 \left( F_{kl} \parbyxx{f_i}{x_k}{x_l} - G_{kl}
 \parbyxx{g_i}{x_k}{x_l} \right) \parby{}{x_i} .\end{multline} All
of them are first and second order operators. Now let's evaluate the
second term
\begin{multline}
 O(0,F) O(0,G) = \ssum{i,j}{1}{n} F_{ij} \parbyxx{}{x_i}{x_j} \biggl(
 \ssum{k,l}{1}{n} G_{kl} \parbyxx{}{x_k}{x_l} \biggr) =\\
 \ssum{i,j,k,l}{1}{n} \biggl( F_{ij} \parbyxx{G_{kl}}{x_i}{x_j}
 \parbyxx{}{x_k}{x_l} + (F_{ij} + F_{ji}) \parby{G_{kl}}{x_j}
 \parbyxxx{}{x_i}{x_k}{x_l} + F_{ij} G_{kl}
 \parbyxxxx{}{x_i}{x_j}{x_k}{x_l} \biggr) .\end{multline}
\begin{multline}
 O(0,F) O(0,G) - O(0,G) O(0,F) = \ssum{i,j,k,l}{1}{n} \biggl( \biggl(
 F_{ij} \parbyxx{G_{kl}}{x_i}{x_j} - G_{ij}
 \parbyxx{F_{kl}}{x_i}{x_j}
 \biggr) \parbyxx{}{x_k}{x_l} +\\
 \biggl( (F_{ij} + F_{ji}) \parby{G_{kl}}{x_j} - (G_{ij} + G_{ji})
 \parby{F_{kl}}{x_j} \biggr) \parbyxxx{}{x_i}{x_k}{x_l} \biggr)
 .\end{multline} Unfortunately the last term
\begin{align}
 \label{eq:60}
 \biggl( (F_{ij} + F_{ji}) \parby{G_{kl}}{x_j} - (G_{ij} + G_{ji})
 \parby{F_{kl}}{x_j} \biggr) \parbyxxx{}{x_i}{x_k}{x_l}
\end{align}
is of third order and, in general, it does not vanish.  Thus we have
shown that the commutator of two general second order operators is of
third order. Consequently, the Leibniz rule simplifications used in
the deterministic case cannot be applied to the general stochastic
linearization problem.

%======== Subsection ================================================
\section{Stochastic Case } 
\label{sub:stoclass}

Since there are two  definitions of coordinate transformation of stochastic
differential equations (\Ito, Stratonovich) and three
definitions of linearity ($g,\sigma,g\sigma$), we face at least six
stochastic problems per a deterministic one. In this section we will
discuss all of them giving at least partial solutions to the feedback
linearization problem. We consider mainly the SISO problem except for
cases where the MIMO extension is trivial.

%======== Subsection ================================================
\subsection{Stratonovich $g$-linearization}
\label{sub:stotech}

We show that the method known for deterministic systems can
be applied without modifications.

\begin{prop} \label{prop:stratgprop}
 The Stratonovich dynamical system~\ssdefmi is $g$-linearizable if and
 only if the deterministic system~\dsdefmi is linearizable. These two
 linearizing transformation{}s are equal. This holds both for SISO and
 MIMO systems and both for SFB and SCT linearization.
\end{prop}
\begin{proof}
 The comparison of transformation laws for deterministic and
 Stratonovich systems shows that the coefficients of $f$ and
 $g$ transform in the same way.  The
 controllability conditions and the definition of linearity are
 also identical (compare Definition~\ref{def:linear} with
 Definition~\ref{def:linearsto}). Identical problems have
 identical  solutions.
\end{proof}

%======== Subsection ================================================
\subsection{Stratonovich $g\sigma$-linearization}
\label{sub:ssgsigma}

The Stratonovich problems are not complicated by the second order \Ito
term. The transformation laws for Stratonovich systems are the same as
the deterministic transformation laws, therefore many results of the
deterministic linearization theory can be used. 

For example, the stochastic SCT $g\sigma$-linearization of a
Stratonovich system is equivalent to the linearization of a
deterministic, non-square MIMO system with two inputs and a single
output. The SFB problem, which is studied in this section, is not as
simple as the SCT one because the feedback influences only the control
input~$u$ (Figure~\ref{fig:picassym}). The ``dispersion input'' is not
a part of the feedback.  Consequently, in order to solve the
Stratonovich SFB $g\sigma$-linearization, we have to deal with a
combined deterministic SFB-SCT problem.
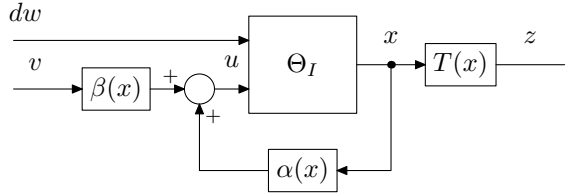
\psfig{fig:picassym}{Asymmetry of  SFB $g\sigma$-Linearization}{picassym}

\subsubsection{Canonical Form}
Recall that we require~$g$-controllability of the resulting system

Since this is a Stratonovich problem, the transformed
vectorfields~$\tilde f$ and~$\tilde g$ do not depend on the dispersion
vectorfield~$\sigma$. Therefore, the {\em control part\/} and the {\em
  dispersion part\/} can be studied independently.

Any~$g$-linear system can be transformed into integrator chain by a
combination of a linear coordinate transformation and linear
feedback. Therefore, if we set~$\sigma=0$, the canonical form is the
integrator chain. 

In general, the dispersion vectorfield~$\tilde \sigma$ is assumed to
be arbitrary constant vectorfield~$\tilde \sigma(x)_i=s_i$, $1\le i\le n$
(See Definition~\ref{def:linearsto}) and this form is preserved by arbitrary linear
transformations. Therefore the canonical form can be written as:
\begin{align}
\label{eq:1000}
\tilde f_i(x) &= x_{i+1} \firstfor 1\le i\le n-1\\
\tilde f_n(x) &= 0\\
\tilde g_i(x) &= 0 \firstfor 1\le i\le n-1\\
\tilde g_n(x) &= 1\\
\label{eq:1077}
\tilde \sigma_i(x) &= s_i \firstfor 1\le i\le n
.\end{align}

We can compare this canonical form with the equations which
define the transformed system~$\tilde \Theta$.

\begin{prop} \label{prop:s1sfb1a}
There is a SFB $g\sigma$-linearizing transformation
 of the SISO Stratonovich system
 \ssdef into a~$g$-controllable linear system if and only if
 there is a solution~$\lambda \RntoR$ of the set of partial differential equations:
\begin{alignat}3
\label{eq:7}
\biglie{\adfg{i}}{\lambda }&=0\firstfor 0 \le i \le n-2\\
\label{eq:13}
\biglie{\adfg{{n-1}}}{\lambda }&\ne 0\\
\label{eq:14}
\biglie{\ad{f}{i}{\sigma}}{\lambda }&=s'_{i+1} \nextfor 0 \le i \le
n-1
\end{alignat}
such that~$s'_i \in \RR$ are constants on~$U$ for $1 \le i \le n$.
Then the linearizing transformation is given by:
\begin{alignat}3
\label{eq:22}
T_i&=\multilie{f}{{i-1}}\lambda  \firstfor 1\le i\le n\\
\alpha &=\frac{-\multilie{f}{{n}}\lambda 
}{\lie{g}{}\multilie{f}{{n-1}}\lambda } &\qquad\qquad&
\label{eq:28}
\beta =\frac{1}{\lie{g}{}\multilie{f}{{n-1}}\lambda }
.\end{alignat}
\end{prop}
\begin{proof}
Assume that~$\Theta_S$ is transformed by~$\combinedtab$ into 
$\tilde \Theta \define \left(\tilde f,\tilde g,\tilde
  \sigma,T(U),T(\origi{x})\right)$ where the~$i$-th components
of~$f$,$g$, and~$\sigma$ can be expressed as:~$\tilde f_i =
\lie{f}{T_i}$, $\tilde g_i =
\lie{g}{T_i}$, $\tilde \sigma_i =
\lie{\sigma}{T_i}$. Moreover, the feedback is defined by~$ u = \alpha  + \beta  v $.
The equations of~$\Theta$ can be compared to
 the equation of the canonical  form~\eqref{eq:1000}-\eqref{eq:1077}.
\begin{alignat}{3}
 \label{eq:160}
 \lie{f}{T_i} &= T_{i+1} \firstfor 1 \le i \le n-1\\
 \label{eq:161}
 \lie{g}{T_i} &= 0 \nextfor 1 \le i \le n-1\\
 \label{eq:162}
 \lie{g}{T_n} &= 1/\beta  \ne 0 \\
 \label{eq:163}
 \lie{f}{T_n} &= -\alpha /\beta  .\end{alignat} 

The
relations~\eqref{eq:7}, \eqref{eq:13}, \eqref{eq:22}, and~\eqref{eq:28} are equivalent to
relations~\eqref{eq:8}-\eqref{eq:12} from 
Proposition~\ref{prop:d1sfb1}. The relation~\eqref{eq:14} can be
verified in a similar way:
\begin{alignat}{3}
 \label{eq:165}
 \lie{\sigma}{T_i} &= s_{i} \firstfor 1 \le i \le n
\end{alignat}
thus by ~\eqref{eq:160}
\begin{alignat}{3}
 \label{eq:254}
 \lie{\sigma}{\lie{f}{T_i}} &= s_{i+1} \firstfor 1 \le i \le n-1
\end{alignat}
 and by ~\eqref{eq:63}
\begin{alignat}{3}
 \label{eq:255}
 \lie{\sigma}{\lie{f}{T_i}} &= \lie{f}{\lie{\sigma}{T_i}}
 -\lie{[f,\sigma]}{T_i} \firstfor 1 \le i \le n-1
\end{alignat}
since the Lie derivative of a constant is zero:
\begin{alignat}{3}
 \label{eq:256} \lie{f}{\lie{\sigma}{T_i}} &= \lie{f}{s_i} = 0\\
 s_{i+1} \define \lie{\sigma}{\lie{f}{T_i}} &=
 -\lie{[f,\sigma]}{T_i} \firstfor 1 \le i \le n-1 .\end{alignat}
The equations~\eqref{eq:14} are obtained by successive application of
this relation. The symbols~$s_i$ are equal to~$s'_{i}$ except for the signs.
\end{proof}

\subsubsection{Conditions for the Control Part}

The necessary conditions for linearizability of the {\em control
  part\/} of~$\Theta$ (\ie the system~$\left(f,g,0,U,\origi{x}
\right)$) can be expressed in geometrical form. We intentionally omit
the dispersion part, using the fact that the resulting system must be
linear when the noise is zero.

Further, the class of all solutions of this subproblem will be
called~$C$. This class can be studied to find  if some member
of~$C$ linearizes the {\em dispersion part\/} of the system.

We can find a geometrical criterion similar to the conditions of
Proposition~\ref{prop:d1sfb2}. In this case these conditions are
necessary but not sufficient since also \eqref{eq:14} must be
satisfied.

\begin{prop}
 \label{prop:s1sfb1}
 SFB~$g\sigma$-linearizing transformation of the Stratonovich system
 $\Theta_S$ into a~$g\sigma$-controllable system 
 linear system exists only if the distribution~$\distrofg{n-2}$ is
 involutive and the distribution~$\distrofg{n-1}$ is $n$-dimensional.
\end{prop}
\begin{proof}
 This theorem is equivalent to Proposition~\ref{prop:d1sfb2} which is
a direct consequence of Proposition~\ref{prop:d1sfb1a} which
corresponds to Proposition~\ref{prop:s1sfb1a}. 
\end{proof}

\subsubsection{Condition for the Dispersion Part}
\label{prop:cftdp}
The conditions of Proposition~\ref{prop:s1sfb1} can be written
in matrix form. We are looking for~$T_1 = \lambda\RntoR$ such that
\begin{align}
 \label{eq:189}
 \left[\,\,\begin{matrix}
 \adfg{0}\\
 \adfg{1}\\
 \vdots\\
 \adfg{n-2}\\
 \end{matrix}\,\,\right]
\quad \left[\,\,\begin{matrix}
 \parby{\lambda }{x_1}\\
 \parby{\lambda }{x_2}\\
 \vdots\\
 \parby{\lambda }{x_n}\\
 \end{matrix}\,\,\right]
= \left[0\right]
.\end{align}

The vectors~$ \adfg{0} \dots \adfg{n-2}$ are written in coordinates
as~$1\times n$ rows. The first matrix is~$(n-1)\times n$. Moreover it is
required that
\begin{equation}
\label{eq:1020}
\biglie{\adfg{n-1}}{\lambda}
\end{equation} is nonzero. 

We will use the algorithm for SFB deterministic linearization (see
Section~\ref{sec:detcase}) to find such a transformation~$\lambda$.Then we will
verify if the conditions for linearity of the dispersion part of the
system~\eqref{eq:14} are also valid. There are~$n$ additional linearity
conditions ($s_i$ are constants):

\begin{align}
\label{eq:1001}
 \left[\,\,\begin{matrix}
 \ad{f}{0}{\sigma}\\
 \ad{f}{1}{\sigma}\\
 \vdots\\
 \ad{f}{n-1}{\sigma}\\
 \end{matrix}\,\,\right]
\quad \left[\,\,\begin{matrix}
 \parby{\lambda }{x_1}\\
 \parby{\lambda }{x_2}\\
 \vdots\\
 \parby{\lambda }{x_n}\\
 \end{matrix}\,\,\right]=
\left[\,\,\begin{matrix}
 s_1\\
 s_2\\
 \vdots\\
 s_n\\
 \end{matrix}\,\,\right]
.\end{align}

In the deterministic case we were satisfied with {\em arbitrary\/}
solution~$\lambda $ to the equations~\eqref{eq:189},
and~\eqref{eq:1020} . In this stochastic case we must find the class
of {\em all\/} solutions and then check if this class contains the
solution for the~$\sigma$ part~\eqref{eq:1001}. Details depend on the
methods used for solving the set of PDEs.

This result is summarized in the following
algorithm:

\begin{itemize}
\item{Find $\Delta_k \define \adfg{i}$ for $0\le i \le k-1$.}
\item{Verify that $\dim(\Delta_{n})$ is $n$.}
\item{Verify that $\Delta_{n-1}$ is
 involutive (see \citet{nijmeijer94} Remark following Definition 2.39)
 otherwise no linearizing 
transformation exists.}
 \item{Find all~$\lambda$ satisfying~\eqref{eq:192} by solving 
PDEs~\eqref{eq:192}; denote~$C$ the set of all
 such functions.} 
\item{Verify that there is a~$\lambda_1\in
 C$ such that the conditions~\eqref{eq:1001} are
 satisfied, otherwise no linearizing
 transformation exists.} 
\item{Compute $T$,$\alpha $,$\beta $ from
 \eqref{eq:160} -- \eqref{eq:163}.}
\end{itemize}

Now, we can illustrate one possible practical approach  which
worked for several simple problems solved by us (see the example in
Section~\ref{sub:appcrane}).

First we can compute the kernel of the matrix~$M_g$ to find the
form~$\omega = \left[ \omega_1, \omega_2, \dots, \omega_n
\right]^T$ which satisfies~$M_g \omega = 0$, \ie~$\omega$ is
perpendicular to~$M_g$. In modern computer algebra systems there is a
single command for this.

Proposition~\ref{prop:s1sfb1} assumes that~$n$
vectorfields~$\Delta_{n} \define \distrofg{n-1}$ form an~$n$
dimensional space. The vectorfields~$\Delta_{n-1} \define
\distrofg{n-2}$ are chosen from them and consequently must form
an~$(n-1)$-dimensional space. Thus their kernel~$d\lambda $ is
exactly one dimensional and arbitrary~$\omega' = c(x)\omega(x)$ also
belongs to the kernel ($c(x)$ is a scalar). 

But not every~$\omega'$ that is perpendicular to~$M_g$ is a
solution to the original linearization problem. The
function~$\omega'$ must be an exact one-form \ie there must be
a scalar function~$\lambda $ such that~$d \lambda =
c(x)\omega(x)$. The Frobenius theorem guarantees that
if~$\Delta_{n-1}$ is involutive, then there is always~$c(x)\in\RR$
such that~$c(x)\omega(x)$ is the exact one-form.

A necessary condition for a one-form~$\omega =
\ssumi{n} \omega_i$ to be exact is
\begin{equation}
 \label{eq:102}
 \parby{\omega_i}{x_j} = \parby{\omega_j}{x_i} \qquad \text{for}
 \qquad 1 \le i,j \le n
.\end{equation}

Hence for every~$1\le i,j \le
n$
\begin{alignat}{3}
 \parby{}{x_j}\left(c(x)\omega_i \right) &=
 \parby{}{x_i}\left(c(x)\omega_j \right),\end{alignat} thus for
every~$1\le i,j \le n$
\begin{alignat}{3}
\label{eq:192}
\parby{c(x)}{x_i}\omega_j - \parby{c(x)}{x_j}\omega_i + c(x)
\left( \parby{\omega_j}{x_i} - \parby{\omega_i}{x_j} \right) = 0
.\end{alignat} The later condition is a set of linear PDEs, with
unknown~$c(x)$, which are guaranteed to have a solution by the
involutivness of~$\Delta_{n-1}$ (the Frobenius theorem). 

In our computations the equation~\eqref{eq:192} was in a simple form
which allowed to determine all the solutions easily. More complicated
cases will require more sophisticated analysis.

%======== Subsection ================================================
\subsection{\Ito $g\sigma$-linearization}
\label{sub:isgsigma}

In the previous subsection we tried to find~$g\sigma$-linearizations
for Stratonovich dynamical systems. Once this is done, the
correcting mapping can be used to
construct \Ito $g\sigma$-linearizing
transformation. This method works for both the
SFB and the
SCT case.

Given an \Ito system~$\Theta_I$, the corresponding
Stratonovich system $\Theta_S$ can be
obtained using the correcting mapping
$\Theta_S = \Corr{\sigma}\left({\Theta_I}\right)$.  Afterward, the Stratonovich
$g\sigma$-linearization algorithm can be applied giving a
linear system~$\Theta_{2S}$. Due to linearity of the
drift
vectorfield~$\tilde \sigma$ of~$\Theta_{2S}$, the
correcting term~$\corr{{\tilde \sigma}}{z}$ of
the backward transformation~${\Corr{{\tilde \sigma}}{}}^{-1}$
vanishes.

\begin{theo} \label{prop:gsigmaprop}
 The SFB 
 $g\sigma$-linearizing
 transformation~$\combinedz_I$ of the~\Ito dynamical system \isdefmi,
 $f(\origi{x})=0$, $\corr{\sigma}{\origi{x}}=0$, into a $g\sigma$-controllable linear
 system  exists if and
 only if there is a SFB
 $g\sigma$-linearizing
 transformation~$\combinedz_S$ of the Stratonovich dynamical
 system

\begin{align}
\label{eq:24}
\Theta_S &= \ssys{\alterbar f}{g}{\sigma}{x} =
 \Corr{\sigma}{(\Theta_I)}\\
\alterbar f &= f + \corr{\sigma}{x}
.\end{align}
Moreover~$\combinedz_I = \combinedz_S \compose \Corr{\sigma}{}$.

\end{theo}
\begin{proof}[Proof (sufficiency)]
 We use the properties of the correcting term
 (Subsection \ref{sub:corr}). Assume that there is a
 mapping~$\combinedz_S$ which transform $\Theta_S$ into a linear
~$g$-controllable system $(Ax,B,S,U,0)$  . By~(\ref{eq:19})
 \begin{align}
 \label{eq:248}
 \combinedz_I = \Corr{{\tilde\sigma}}{}^{-1} \compose \combinedz_S
 \compose \Corr{\sigma}{} 
.\end{align}
 The backward correcting
transformation $\Corr{{\tilde \sigma}}^{-1}$ is identity because the correcting term
of a linear mapping $\corr{{\tilde \sigma}}{x}$ is zero. Thus~$
\combinedz_I =  \combinedz_S \compose \Corr{\sigma}{} $
and ~$\combinedz_I (\Theta_{I})$ equals $(Ax,B,S,U,0)$, which is
linear and~$g$-controllable by assumption.
\end{proof}
\begin{proof}[Proof (necessity)]
  Assume that there is the \Ito transformation~$\combinedz_I$ which
  linearizes $\Theta_I$ and by~(\ref{eq:24}) $\Theta_I =
  {\Corr{\sigma}{}}^{-1}(\Theta_S)$. Construct Stratonovich
  linearization by~$\combinedz_S = \combinedz_I \compose
  {\Corr{\sigma}{}}^{-1} $.  Hence~$\combinedz_I$ linearizes
  ${\Corr{\sigma}{}}^{-1}(\Theta_S)$ and~$\combinedz_S$
  linearizes~$\Theta_S$ into the same linear and controllable system
  as~$\combinedz_I$.
\end{proof}

%======== Subsection ================================================
\subsection{\Ito~$g$-linearization}
\label{sub:sfbig}

The \Ito $g$-linearization problem is probably the most complicated
variant of exact linearization studied in this paper. The dispersion
vectorfield of an \Ito dynamical system transformed by a coordinate
transformation~$\sctt$ consists of two terms: the transformed
vectorfield~$\tantra{T} f$ and the \Ito term~$\ito{\sigma}$. We
require that the sum of these  terms is linear, thus the
nonlinearity of the drift~$\tantra{T}f$ must compensate for the \Ito
term. Since the \Ito term behaves to~$T$ as a second order
differential operator, this problems generates a set of second order
partial differential equations. One can attempt to use simplifications
as in the deterministic linearization, namely, the recursive Leibniz
rule~\eqref{eq:10}. Unfortunately, this approach does not work for the
stochastic case. In general, the \Ito equations cannot be easily
simplified by commutators because the commutator of second order
operators is not a second order operator but a third order operator
(see Subsection~\ref{sub:stoinvar}).

Nevertheless, there are special cases for which simpler conditions
 can be found. The most important special case (commuting $g$
 and $\sigma$) will be studied here.

\subsubsection{Canonical Form ---$n$ unknowns}

The canonical form for the $g$-linearization
  is the integrator chain with a nonlinear drift
\begin{align}
\label{eq:1007}
\tilde f_i(x) &= x_{i+1} \firstfor 1\le i\le n-1\\
\label{eq:1008}
\tilde f_n(x) &= 0\\
\label{eq:1009}
\tilde g_i(x) &= 0 \firstfor 1\le i\le n-1\\
\label{eq:1010}
\tilde g_n(x) &= 1
.\end{align}
Assume that there is a~$g$-linear
system~$\Theta_I = (Ax,B,\sigma(x),U,\origi{x})$. Then the drift part
of~$\Theta_I$ can be transformed by a {\em linear} transformation into
the integrator chain. This is because the \Ito term of a linear
transformation vanishes. 

The equations which define~$T$ can be obtained by comparing  this
canonical form with the equations of~$\tilde \Theta$.

\begin{prop} \label{prop:s1sfb2}
  Let~\isdef be an \Ito dynamical system with~$f(\origi{x})=0$ such
  that~$\corr{\sigma}{\origi{x}}=0$. There is a SFB~$g$-linearizing
  transformation~$\combinedtab$ of the system~$\Theta_{I}$ into a
  ~$g$-controllable linear system if and only if there is a solution
  $T_i\RntoR$, $1\le i \le n$, to the set of partial differential
  equations defined on~$U$:
\begin{alignat}{3}
\label{eq:50}
T_{i+1}&=\lie{f}{T_i}+\ito{\sigma}{T_i} \firstfor 1\le i\le{n-1}\\
\label{eq:51}
\lie{g}{T_i}&=0 \nextfor 1\le i\le{n-1}\\
\label{eq:52}
\lie{g}{T_n}&\not=0 .\end{alignat} The symbol~$\ito{\sigma}$
denotes the \Ito operator (see Definition~\ref{eq:16}). The
feedback can be constructed as:
\begin{align}
\label{eq:53}
\alpha &=-\frac{(\lie{f}{T_n}+\ito{\sigma}{T_n})}{\lie{g}{}T_n}
\qquad\qquad \beta =\frac{1}{\lie{g}{T_n}}
.\end{align}
\end{prop}
\begin{proof}

The~$i$-th components
of~$f$,$g$, and~$\sigma$ are:~$\tilde f_i =
\lie{f}{T_i} + \ito{\sigma}{T_i}$, $\tilde g_i =
\lie{g}{T_i}$, $\tilde \sigma_i =
\lie{\sigma}{T_i}$. 
 The partial differential equations \eqref{eq:50}, \eqref{eq:51}
 and \eqref{eq:52} are obtained by comparison of \eqref{eq:86}-\eqref{eq:89}
 with the equations~\eqref{eq:1007}-\eqref{eq:1010}.
\end{proof}
\subsubsection{PDEs of single unknown}

One can attempt to reduce the equations~\eqref{eq:50},
\eqref{eq:51}, and~\eqref{eq:52}, to a set of equations of a
single unknown, similarly to the results of Proposition \ref{prop:d1sfb1}.

\begin{coroll}
 \label{prop:sfbooo} Define the general second order
 operator~$O(f,F)$ as in~\eqref{eq:17}. The exponential notation
 for~$O$ will be defined recursively: $O^0(f,F)T \define T$ and
 $O^{l+1}(f,F)T \define O(f,F)O^l(f,F)T$ for~$l \ge 0$. Next, define
 \begin{align} F_{ij} \define \frahalf \sigma_i \sigma_j .
 \end{align} Then the set of partial differential equations~\eqref{eq:50}-\eqref{eq:52} of~$n$ unknowns has a 
 solution if and only if there is a 
 solution~$\lambda \RntoR$ defined on~$U$ to the set of
 PDEs of single unknown:
 \begin{alignat}{3} \label{eq:153}
 O(g,0)O^i(f,F)\lambda &=0\firstfor 0\le i\le n-2\\ 
 \label{eq:154} O(g,0)O^{n-1}(f,F)\lambda &\ne 0 . \end{alignat}
 The original solution and the feedback
 can be found as \begin{alignat}{3}
 T_i&=O^{i-1}(f,F)\lambda \firstfor 1 \le i \le n\\
 \label{eq:58} \alpha &= -
 \frac{O(f,F)^{n-1}\lambda }{O(g,0)^{n-1}\lambda } \quad\quad \beta =
 \frac{1}{O(g,0)^{n-1}\lambda } . \end{alignat}
\end{coroll}
\begin{proof}
 Since~$T_1=\lambda $ and by definition
 of~$\ito{\sigma}{}$ \eqref{eq:16} and~$O$~\eqref{eq:17}:
 \begin{alignat}3
\label{eq:264}
O(f,F)T_i &= \lie{f}{T_i} + \ito{\sigma}{T_i} \firstfor 1\le i\le n\\
\label{eq:265}
O(g,0)T_i &=\lie{g}{T_i} \nextfor 1\le i\le n ;
 \end{alignat} then~$T_{i+1} = O(f,F)T_i$ by~\eqref{eq:50} and
 $T_i = O^{i-1}(f,F)T_1 = O^{i-1}(f,F)\lambda $. Similarly
 the equation \eqref{eq:53} implies~\refprop{eq:58}.
\end{proof}

Note, that for the deterministic case~$\sigma=0$, the
operators~$O(g,0)$ and~$O(g,0)O^i(f,F)$ degenerate
to~$\lie{g}{}$ and to~$\lie{g}{}\multilie{f}{i}{}$ respectively,
thus the result is the same as that of Proposition~\ref{prop:d1sfb1}.

In general, the equations of the system are of an order up to~$2n$ and
cannot be reduced to a lower order. The
commutator of two second order operators is of third
order as was pointed out by~\eqref{eq:60}. In particular, for~$i=1$ we
have to evaluate~$C(g,0,f,F) \definer O(a,A)$, which {\em is\/} of
second order due to the fact that~$G=0$. But starting from~$i=2$ the
commutator~$C(f,F,a,A) = C(f,F,C(g,0,f,F))$ is of
third order.

\subsubsection{Correcting Term}

The same problem can be reformulated using conversion to the
Stratonovich formalism. One can compute the drift vectorfield~$\alterbar f$
of the equivalent Stratonovich system by applying the
correcting term~$\alterbar f \define f +
\corr{\sigma}{x}$. Then the Stratonovich system~$\ssys{\alterbar
 f}{g}{\sigma}{x}$ may be transformed, by a suitable transformation,
to such a form~$(\tilde f,\tilde g,\tilde \sigma,T(U),0)$ that after
applying the backward correcting
term~$-\corr{{\tilde \sigma}}{z}$ the resulting \Ito system will be
linear.

Compare this formulation with the~$g\sigma$-linearization where the
backward correcting term vanished due to
linearity of~$\tilde \sigma$. This does not happen with
the~$g$-linearization, and the backward correcting term is a part of the equations.

In general, it may be difficult to solve these equations. Nevertheless
there are special cases when the solution can be obtained. See for example
Section~\ref{sub:appcrane}. Another important case (commuting $g$~
and~$\sigma$) is studied below.

\begin{coroll} \label{prop:s1sfb3}
 The equations of Proposition~\ref{prop:s1sfb2} are equivalent to~$n$
 partial differential equations:
\begin{alignat}3
\label{eq:143}
T_{i+1}&=\lie{\alterbar f}{T_i} -\corr{{\alterbar \sigma}}{z}=\\
&\lie{\alterbar f}{T_i} + \frahalf \lie{\sigma}{\lie{\sigma}{T_i}} \firstfor 1 \le i \le{n-1}\\
\label{eq:144}
\lie{g}{T_i}&=0\nextfor 1 \le i \le{n-1}\\
\label{eq:145}
\lie{g}{T_n}&\not=0 .\end{alignat} Then the
feedback can be constructed as
\begin{align}
  \alpha &=-\frac{\lie{\alterbar f}{T_n} + \frahalf
    \lie{\sigma}{\lie{\sigma}{T_n}}}{\lie{g}{T_n}} \qquad\qquad
  \beta =\frac{1}{\lie{g}{T_n}} .\end{align} Where $\alterbar f
\define f + \corr{\sigma}{x}$.
\end{coroll}
\begin{proof}
 The equations~\eqref{eq:143} can be obtained
 from~\eqref{eq:50} by applying~\eqref{eq:215} 
 and~\eqref{eq:140}:
\begin{multline}
 \label{eq:216}
 T_{i+1} = \lie{f}{T_i}+\ito{\sigma}{T_i} = \lie{\left(
 \alterbar f - \corr{{\alterbar \sigma}}{x} \right)}{T_i}+
 \ito{\sigma}{T_i} = \\
 \lie{\alterbar f}{T_i} - \lie{\corr{{\alterbar
 \sigma}}{x}}{T_i} + \ito{\sigma}{T_i} = \lie{\alterbar
 f}{T_i} + \frahalf\lie{\sigma}{\lie{\sigma}{T_i}}
 .\end{multline}

The other equations are adopted from Proposition~\ref{prop:s1sfb2}.
The set of PDEs of~$n$ unknowns can be transformed into a set
of PDEs of a single unknown~$\lambda =T_1$, but
the order of the equation will be~$2n-1$.
\end{proof}

\begin{rem} \label{prop:s1sfb4}
  Observe that the set of~$n$ second order partial differential
  equations~\eqref{eq:143}-~\eqref{eq:145} defined in
  Proposition~\ref{prop:s1sfb3} can be transformed, by introducing 
  new variables~$S_i = \lie{\sigma}{T_i}$,  to the following system
  of~$2n-1$ first order partial differential equations for~$1 \le i
  \le {n-1}$:
\begin{alignat}3
\label{eq:146}
\lie{g}{T_i}&=0\\
\label{eq:147}
\lie{\sigma}{T_i}-S_i&=0\\
\label{eq:148}
\lie{{\alterbar f}}{T_i} + \frahalf\lie{\sigma}{S_i}&=T_{i+1}\\
\label{eq:149}
\lie{g}{T_n}&\not=0 .\end{alignat} $T_i$ 
and~$S_i\RntoR$ are unknown real valued functions defined on $U$.
\end{rem}

\subsubsection{Systems with Commuting $g$ and $\sigma$}
There is a special case of \Ito dynamical systems for which the 
solution is completely known and can be computed using
only first order PDEs.

\begin{theo}
 \label{prop:ssfbco}
 Let~\isdef be a SISO \Ito dynamical system. If the
 vectorfield~$\sigma$ commutes with all vectorfields~$\ad{{\alterbar
 f}}{i}{g}$ for~$0\le i\le n-1$, \ie ~$[\ad{{\alterbar
 f}}{i}{g},\sigma]=0$, where $\alterbar f
\define f + \corr{\sigma}{x}$, then the \Ito system is $g$-linearizable if and
 only if the distribution
\begin{align}
\alterbar \Delta_{n} \define
 \sspan\left\{{\ad{{\alterbar f}} {i}{g}, \iva{n-1}}\right\}
\end{align}
is
 nonsingular on~$U$ and the distribution
\begin{align}
\alterbar \Delta_{n-1}
 \define \sspan\left\{{\ad{{\alterbar f}}{i}{g}, \iva{n-2}}\right\}
\end{align}
is involutive on~$U$. If these conditions hold, then a
 solution~$\lambda \RntoR$ to the set of partial differential
 equations exists
\begin{align}
\label{eq:77}
\biglie{\ad{{\alterbar f}}{i}{g}}{\lambda }&=0\rfor{\iva{n-2}}\\
\label{eq:78}
\biglie{\ad{{\alterbar f}}{{n-1}}{g}}{\lambda }&\not=0 .\end{align}
the linearizing transformation is given by:
\begin{alignat}3
\label{eq:79}
T_1&=\lambda \\
\label{eq:80}
T_{i+1}&= \lie{\alterbar f}{T_i} + \frahalf
\lie{\sigma}{\lie{\sigma}{T_i}} \firstfor 1
\le i \le{n-1}\\
\alpha &=\frac{-\multilie{\alterbar f}{{n}}\lambda 
}{\lie{g}{}\multilie{\alterbar f}{{n-1}}\lambda } &\qquad\qquad&
\label{eq:81}
\beta =\frac{1}{\lie{g}{}\multilie{\alterbar f}{{n-1}}\lambda }
.\end{alignat}
\end{theo}
\begin{proof}
 We will apply the Leibniz rule to the relation of~(\ref{eq:143})
Corollary~\ref{prop:s1sfb3} to expand the
 term~$\lie{g}{T_{i+1}}$ for~$1\le i\le n-1$ 
\begin{multline}
 \label{eq:195}
 \lie{g}{T_{i+1}} = \lie{g}{ \left( \lie{{\alterbar f}}{T_i}+
 \frahalf \lie{\sigma}{\lie{\sigma}{T_i}} \right) } =
 \lie{g}{\lie{{\alterbar f}}{T_i}} + \frahalf
 \lie{g}{\lie{\sigma}{\lie{\sigma}{T_i}}}
 =\\
 \lie{{\alterbar f}}{\lie{g}{T_i}} - \lie{[\alterbar f,g]}{T_i}
 + \frahalf \lie{\sigma}{\lie{g}{\left(\lie{\sigma}{T_i}\right)}}
 - \frahalf \lie{[\sigma,g]}{\left(\lie{\sigma}{T_i}\right)}
 =\\
 0 - \lie{[\alterbar f,g]}{T_i} + \frahalf \lie{\sigma}{ \left(
 \lie{\sigma}{\lie{g}{T_i}} - \lie{[\sigma,g]}{T_i} \right)
 } - \frahalf \left( \lie{\sigma}{\lie{[\sigma,g]}{T_i}} -
 \lie{[\sigma,[\sigma,g]]}{T_i} \right)
 =\\
 -\lie{[\alterbar f,g]}{T_i} + 0 -  \frahalf
 \lie{\sigma}{\lie{[\sigma,g]}{T_i}} + \frahalf
 \lie{[\sigma,[\sigma,g]]}{T_i} .\end{multline} If the
vectorfields~$g$ and~$\sigma$ commute, then the second and
third terms vanish. If, moreover, $[\sigma,[\alterbar f,g]] = 0$ then
\begin{alignat}{3}
 \label{eq:257}
 \lie{g}{{T_{i+2}}} &= - \lie{{[\alterbar f,g]}}{ \left(
 \lie{{\alterbar f}}{T_i} + \frahalf
 \lie{\sigma}{\lie{\sigma}{T_i}} \right)} = \lie{{[\alterbar
 f,[\alterbar f,g]]}}{T_i} .\end{alignat} In general
if~$[\sigma,\ad{{\alterbar f}}{i}{g}]=0$ for~$0\le i\le n-1$ then 
\begin{alignat}{3}
 \label{eq:258}
 \lie{g}{T_{k}} &= (-1)^k \lie{{[\ad{{\alterbar f}}{k}{g},]}}{T_1}
 .\end{alignat}

Thus the equations~\eqref{eq:143} and \eqref{eq:144} will be
equivalent to
\begin{alignat}{3}
 \label{eq:193}
 \biglie{\ad{{\alterbar f}}{i}{g}}{\lambda }&=0\firstfor{\iva{n-2}}\\
 \label{eq:194}
 \biglie{\ad{{\alterbar f}}{{n-1}}{g}}{\lambda }&\not=0
,\end{alignat} which are of the same form as the equations of
Proposition~\ref{prop:d1sfb1} and consequently the conditions
from Proposition \ref{prop:d1sfb2} can be used.
\end{proof}

%======== Subsection ================================================
\subsection{\Ito and Stratonovich $\sigma$-linearization}
\label{sub:sfbchsigma}

The stochastic SFB $\sigma$-linearization problem is similar to
deterministic SCT linearization.
The dispersion
vectorfield~$\sigma$ transforms in the same way as deterministic
drift vectorfields do. Consequently, no \Ito term
complicates the transformation. Moreover, the \Ito and
Stratonovich cases are equivalent.

On the other hand, in the SFB $\sigma$-linearization we are free to
choose the feedback~$\feedbackab$ that perturbs the
drift
vectorfield~$f$ into~$\tilde f = f + g \alpha $.

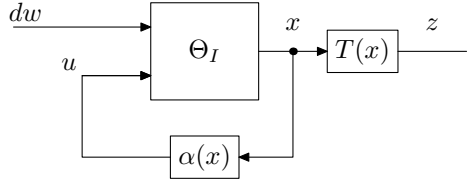
\psfig{fig:sigmasfb}{Ito and Stratonovich
 $\sigma$-linearization}{sigmasfb}
\begin{prop} \label{prop:sfsigma}
 Let~$\Theta$ be a SISO stochastic
 system~$\Theta = \ssysfgx$. There is a
 SFB $\sigma$-linearizing
 transformation~$\combinedtab$ into a $\sigma$-controllable linear
 system if and only if there is a smooth
 feedback function~$\alpha \RntoR$ such that the
 deterministic system~$\dsys{f+g\alpha }{\sigma}{x}$ has a
 SCT
 linearizing
 transformation~$\sctt$. Equivalently,  there must be
such~$\alpha $ that the the modified odd bracket condition:
\begin{equation}
 \label{eq:59}
 [\sigma,\ad{f+g \alpha  }{l}{\sigma}]=0 \qquad\text{for}\qquad l = 1,\dots,2n-1 
\end{equation}
is satisfied (see~\eqref{eq:45}). The resulting combined
transformation consists of composition of the coordinate transformation~$\sctt$ and the
feedback~$\feedbackab$ where~$\beta $ is arbitrary
function of~$x$; for instance~$\beta  = 1$.
\end{prop}
\begin{proof}
Compare definition of linearity of a deterministic system with
 definition of 
 $\sigma$-linearity. The system is is $\sigma$-linearizable if and only if the
 deterministic systems~$\dsys{f+g\alpha }{\sigma}{x}$ is
 SCT linearizable (see Corollary~\ref{prop:s1ctt2a}).
 
 It is evident that the function~$\beta $ (see
 Figure~\ref{fig:sigmasfb}) has no effect on the dispersion part and
 can be chosen arbitrarily. (Probably nonzero for otherwise the system will
 be $g$-uncontrollable).
\end{proof}

The condition~\eqref{eq:59} can be expressed in terms of
derivatives of $\alpha $ using bracket
relations known from differential geometry. For example, for~$l=1$:
\begin{align}
\label{eq:156}
[\sigma,[f+g \alpha ,\sigma ]] &= [\sigma, [f,\sigma ]+[g \alpha
,\sigma ]] = [\sigma,[f,\sigma ]]+[\sigma,[g \alpha ,\sigma ]]=\\
&=[\sigma,g+[f+g,\sigma ]] + g \lie{\sigma}{\lie{\sigma}{\alpha }}
.\end{align}
The other conditions for~$k=3,5,7,\dots$ can be expressed in a similar
way giving the set of~$n$ partial
differential equations of the order up to~$2n$ for example by a
computer using symbolic algebra tools. The
problem is not very interesting from the practical point of view.

%======== Subsection ================================================
\section{Example---Crane}
\label{sub:appcrane}

In this section the methods of stochastic exact linearization are
demonstrated on an example --- control of a crane under
the influence of random disturbances. The description of the plant was
adopted from~\citet{ackermann93} where the model of a crane linearized
by approximative methods was studied. Unlike Ackermann, we control the
same system using the exact model. Moreover the influence of random
disturbances is added. \psfig{fig:crane1}{Crane}{crane1} Consider the
crane of Figure~\ref{fig:crane1}, which can be used for example for
loading containers into a ship. The hook must
be automatically placed to a given position.
Feedback control is needed in order to dampen
the motion before the hook is lowered into the ship. The input signal
is the force~$u$ that accelerates the crab. The crab mass is~$m_C$,
the mass of the load~$m_L$, the rope length is~$l$, and the
gravity acceleration~$g$.

We assume that the driving motor has no nonlinearities, there is no
friction or slip, no elasticity of the rope and no damping
of the pendulum (\eg from air drag). We will define four state
variables: the rope angle $x_1$ (in radian), the angular velocity~$x_2
= \dot x_1$, the position of the crab~$x_3$, and the velocity of
the crab~$x_4 = \dot x_3$. As shown in~\citez{ackermann93}, the
plant is described by two second order differential equations:
\begin{align}
 \label{eq:72}
 u &= (m_L + m_C) \ddot x_3 + m_L l ( \ddot x_1 \cos x_1 - \dot x_1^2
 \sin
 x_1) \\
 \label{eq:73}
 0 &=m_L \ddot x_3 \cos x_1 + m_L l \ddot x_1 + m_L g \sin x_1
.\end{align}
Additionally, we assume that the load is under influence of random
disturbance, which can be modeled as a white noise
process. The disturbance (wind) is horizontal, has zero mean and
 can be described by the \Ito differential~$dw$:
\begin{align}
 \label{eq:74}
 dx_2 = \frac{F \cos x_1}{m_L l}\,dw,\end{align} where~$F$ is a
constant having the physical unit of force.

We used symbolic algebraic system \Mathematica to handle the
computations. The complete \Mathematica worksheet can be downloaded
from the web page of the author  \hyref{http://www.tenzor.cz/sladecek}.

\Mathematica was used to solve the equations of the system for
unknown values~$\dot x_2$ and~$\dot x_4$ (angular and
positional acceleration). Values of vectorfields~$f$, $g$, and
$\sigma$ were derived as follows:

\begin{align}
 \label{eq:76}
 f &= \left[ x_2, -\frac{\sin x_1 \left( g (m_L+m_C) + l m_L x_2 \cos
 x_1 \right) }
 {l (m_C + m_L - m_L \cos^2 x_1) }, x_4, 0\right]^T\\
 g &= \left[ 0, -\frac{\cos x_1 }
 {l (m_C + m_L - m_L \cos^2 x_1) }, 0, u\right]^T\\
 \sigma &= \left[0, \frac{F \cos x_1}{m_L l},0,0 \right]^T
.\end{align}

The state space model is shown in Figure~\ref{fig:crane2}. We can see
that the positional state variables~$x_3$ and~$x_4$ are isolated from
the angular state variables~$x_1$ and~$x_2$. Later, we will
concentrate on the angular variables pretending that the load will be
stabilized no matter where the crane is. Consequently, we obtain only
two-dimensional system on which the exact linearization techniques can
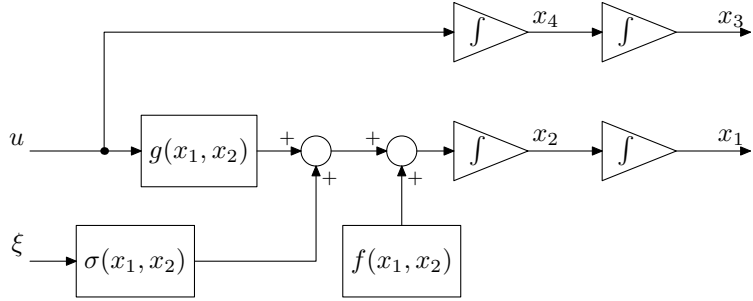
be demonstrated. \psfig{fig:crane2}{The State Space Model of
  Crane}{crane2}

Next, consider the random disturbances. Because the
correcting term~$\corr{\sigma}{x}$ is zero, 
there is no difference in using either the
\Ito or the Stratonovich integral. In case of more ``nonlinear'' noise, one of
the integrals must be selected. If the \Ito model is chosen,
Theorem~\ref{prop:gsigmaprop} must be applied.

Now we evaluate the conditions of Proposition~\ref{prop:s1sfb1}
to check that the system is linearizable. In fact, we must only
evaluate the non-singularity condition because every one-dimensional
distribution is involutive, and the
integrability is satisfied automatically. To this
end, we will compute the null space (kernel) of the
matrix~$[[f,g],g]$, which is empty and therefore the matrix is
nonsingular. We conclude, that the {\em deterministic\/}
SFB problem is solvable. 

Notice, that the system is already in the integrator chain form and
hence~$\lambda = x_1$ satisfies this condition. Therefore, the
{\em deterministic\/} system is linearizable by feedback only, with no
state space transformation at all, \ie $z=T(x)=x$.

This choice of the output function~$\lambda $ is natural but does
not cancel the nonlinearity in the dispersion coefficient~$\sigma$.
For this purpose, we must use the algorithm of Section~\ref{prop:cftdp} to construct
another nontrivial coordinate transformation~$T$.

To obtain this transformation, we must find the space of all
functions~$\lambda $ satisfying conditions for
feedback linearity~\eqref{eq:7}. Observe
that~$\lie{g}{\lambda }$ must be zero hence
\begin{align}
 \label{eq:75}
 \parby{\lambda }{x_1} g_1 + \parby{\lambda }{x_2} g_2 &= 0
 .\end{align}

Since~$g_1=0$ and~$g_2 \ne 0$ in
neighborhood of~$x_0$, then
$\parby{\lambda }{x_2} = 0$ and $\lambda  = c_1(x_1)$ is a
function of~$x_1$ only (\ie without~$x_2$). The coordinate transformation is $T =
\left[\lambda ,\lie{f}{\lambda } \right]^T$. We want to select
such $c_1(x_1)$ that the dispersion
vectorfield~$\tilde \sigma \define \tantra{T}\sigma$ in the new
coordinate system~$z=T(x)$ will be constant:
\begin{align}
 \label{eq:70}
 \parby{c_1(x_1)}{x_1} \frac{F \cos x_1}{m_L l} = \text{constant}
 .\end{align}

We decided to define the constant as~${F}/{(m_L l)}$, therefore
\begin{align}
 \label{eq:221}
 \parby{c_1}{x_1} &= \frac{1}{\cos x_1}
\end{align}
 and
\begin{multline}
 \label{eq:71}
 T_1 = \lambda  = c_1(x_1) = \int \frac{1}{\cos x_1} \, d x_1 =\\
 -\ln \left( \cos \frahalf{x_1} - \sin \frahalf{x_1} \right) +\ln
 \left( \cos \frahalf{x_1} + \sin \frahalf{x_1} \right)
 .\end{multline}
\begin{align}
 \label{eq:222}
 T_2 = \lie{f}{\lambda } &= x_2 \sec x_1
.\end{align}

Finally, we can compute the feedback from~\eqref{eq:28}. In the
\Mathematica worksheet we validate the results by computing~$\tilde \Theta
= \combinedtab \ssysfgx $, The computation showed that the
system~$\ssys{\hat f}{\hat g}{\hat \sigma}{x}$ is in the integrator
chain form in the~$z$ coordinate chart.

\begin{multline}
b=\frac{1}{l \left( m_c + m_l \left( \sin (x_1) \right)^2 \right)}\\
a=\tan ({x_1}) \left( \sec (x_1)\,x_2^2 -
         b\,g\, ( m_c +  m_l ) + l\,m_l\,{x_2}\cos ({x_1}) \right)
.\end{multline}

\section{Conclusion}

\subsection{Main Results}

\begin{romanenu} 
 
\item The structure of the stochastic linearization problem is much richer
 than the structure  of the deterministic one. Two definitions of coordinate
 transformation exist and there are differences between the~$g$,~$\sigma$,
 and~$g\sigma$-linearization.

\item In the case of \Ito integrals, the coordinate transformation
 laws are of second order (the \Ito rule). 

There is a large difference between
$g\sigma$-linearization and $g$-linearization. In the former case the
effect of the \Ito term can be reduced to the first order operator and
consequently the problem is solvable by differential geometry. On the
other hand, in the later case there is no easy method to elliminate
the \Ito term and the a set of second order partial differential
equations must be solved to get the linearizing transformation.

\item We have given (at least partial) solutions to all
 SISO SFB problems. The results are listed in Table~\ref{tab:sfb}.

\begin{table}[htpb]
 \begin{center}
\begin{tabular}{llll}
\hline
Linearization&$g$&$g\sigma$&$\sigma$\\
\hline
Deterministic&\reftpropag{prop:d1sfb2}&&\\ 
Stratonovich &\refppropag{prop:stratgprop}
&\reftpropag{prop:s1sfb1} 
&\refppropag{prop:sfsigma}\\ 
\Ito&\refcorol{prop:sfbooo}
&\reftpropag{prop:gsigmaprop}
&\refppropag{prop:sfsigma}\\ 
\hline

\end{tabular}
\caption{Overview of results --- SISO SFB case}
\label{tab:sfb}
\end{center}
\end{table}

\item \Ito linearization problems can be approached by means of the
 correcting term. The \Ito differential
 equation can be converted to the Stratonovich equation
 whose behavior under coordinate transformations is simpler.
 This method is only partially applicable to the $g$-linearization.
\item An important special case was identified for the \Ito
  $g$-linearization. The case is characterized by commuting control
  vectorfields~$g$ and dispersion
  vectorfields~$\sigma$.  Solutions can be found
  using  first order methods.
\item Computer algebra proved to be a useful tool for solving exact
 linearization problems. 
\item Industrial applications of the exact linearization in general
 are still unlikely, mainly due to complexity, sensitivity, and
 limited robustness of the control laws designed by the method. 
\end{romanenu}

\subsection{Future Research}

\begin{romanenu}
\item Find a solution to the \Ito~$g$-linearization problem in general
 case, including geometric criteria, using second order geometry.
\item Analyze the computability issues; implement a universal symbolic
  algebra toolbox for the problem.
\item Solve the SCT problem.
\item Extend the results to the MIMO systems.
\item Extend the results to the input-output problems and
  linearization of autonomous systems. Work out the applications of
  nonlinear filtering.
\item Perhaps, some of the results can be used as a starting
 point for approaching more general class of problems as the problems
 of disturbance decoupling, input invariance of stochastic
 non--linear systems, or problems of reachability and
 observability.

\end{romanenu}

\bibliography{sladecek}
\bibliographystyle{plainnat}

\end{document}